\documentclass[11pt]{amsart}
\usepackage{amssymb,amsmath,amsthm}
\usepackage{graphicx,diagrams}

\usepackage{a4wide}

\numberwithin{equation}{subsection}

\newtheorem{thm}{Theorem}[subsection]
\newtheorem{prop}[thm]{Proposition}
\newtheorem{lem}[thm]{Lemma}
\newtheorem{cor}[thm]{Corollary}
\newenvironment{rem}{\vspace{3mm}\noindent
{\bf Remark.}}{\vspace{3mm}}

\newcommand{\Pf}{\noindent {\it Proof}}

\newcommand*{\Dd}{\mathop{\mathrm D\kern0pt}\nolimits}

\long\def\comment#1{}

\newcommand{\Tr}{\operatorname{Tr}}

\newcommand{\OO}{{\mathcal O}}

\newcommand{\DD}{{\mathcal D}}

\renewcommand{\Re}{\operatorname{Re}}

\newcommand{\lan}{\langle}
\newcommand{\ran}{\rangle}
\newcommand{\Coh}{\operatorname{Coh}}

\newcommand{\de}{\delta}
\newcommand{\eps}{\epsilon}
\renewcommand{\ker}{\operatorname{ker}}

\newcommand{\id}{\operatorname{id}}

\newcommand{\ov}{\overline}
\newcommand{\we}{\wedge}
\renewcommand{\Im}{\operatorname{Im}}

\newcommand{\EE}{{\mathcal E}}

\newcommand{\Om}{\Omega}
\newcommand{\dbar}{\overline{\partial}}
\newcommand{\pa}{\partial}
\newcommand{\Hom}{\operatorname{Hom}}

\newcommand{\Ext}{\operatorname{Ext}}

\renewcommand{\a}{\alpha}
\renewcommand{\b}{\beta}

\newcommand{\la}{\lambda}
\renewcommand{\th}{\theta}
\newcommand{\C}{{\Bbb C}}

\newcommand{\R}{{\Bbb R}}
\newcommand{\Z}{{\Bbb Z}}

\newcommand{\La}{\Lambda}

\newcommand{\wt}{\widetilde}
\newcommand{\ot}{\otimes}

\newcommand{\sub}{\subset}
\newcommand{\ed}{\qed\vspace{3mm}}

\begin{document}

\title{$A_{\infty}$-algebras associated with elliptic curves and Eisenstein-Kronecker series}
\author{Alexander Polishchuk}
\address{Department of Mathematics, University of Oregon, Eugene, OR 97403}  
\email{apolish@uoregon.edu}
\thanks{Supported in part by the NSF grant DMS-1400390}

%\date{}

\begin{abstract} We compute the $A_{\infty}$-structure on the self-$\Ext$ algebra of the vector bundle $G$ over an elliptic curve
of the form 
$G=\bigoplus_{i=1}^r P_i\oplus \bigoplus_{j=1}^s L_j$, where $(P_i)$ and $(L_j)$ are line bundles of degrees $0$ and $1$,
respectively. The answer is given in terms of Eisenstein-Kronecker numbers $(e^*_{a,b}(z,w))$. The $A_\infty$-constraints
lead to quadratic polynomial identities between these numbers, allowing to express them in terms of few ones. 
Another byproduct of
the calculation is the new representation for $e^*_{a,b}(z,w)$ by rapidly converging series.  
\end{abstract}

\maketitle

\section*{Introduction}

Let $\La\sub \C$ be a lattice. Let $a(L)$ denote the area of the fundamental parallelogram of $\La$, and set
$A(\La)=\frac{a(\La)}{\pi}$.
Recall (see \cite[ch.\ VIII]{Weil}) that the {\it Eisenstein-Kronecker-Lerch series} are given by
$$K^*_a(z,w,s;\La)=\sum_{\la\in\La\setminus\{-z\}}\frac{(\ov{z}+\ov{\la})^a}{|z+\la|^{2s}}\lan \la,w\ran_\La,$$
where $a\in\Z_{\ge 0}$, $z,w\in\C\setminus\La$, $s$ is a real number,
$$\lan z,w\ran_\La:=\exp[(z\ov{w}-w\ov{z})/A(\La)].$$
This series converges absolutely for $\Re s>a/2+1$. It is known that $K^*_a(z,w,s;\La)$
analytically extends (for fixed $z,w$) to a meromorphic function on the entire $s$-plane, with possible poles
only at $s=0$ (for $a=0$, $z\in \La$) and at $s=1$ (for $a=0$, $w\in \La$). 
Using this analytical continuation  
the {\it Eisenstein-Kronecker numbers} $e^*_{a,b}(z,w;\La)$, for integers $a\ge 0$, $b>0$, 
are defined as the following special values:
$$e^*_{a,b}(z,w)=K^*_{a+b}(z,w,b;\La).$$
Note that $e^*_{a,b}(z,w)$ descends to a section of a line bundle on $C\times C$, where $C=\C/\La$ is the corresponding
elliptic curve. This section is not a continuous, however, it has continuous (even real analytic) restrictions to all the strata
of the stratification of $C\times C$ induced by the stratification $C=\{0\}\sqcup (C\setminus\{0\})$.

We refer to \cite[ch.\ VIII]{Weil} and \cite{Bannai} for basic properties of Eisenstein-Kronecker numbers.
Note that they are of number-theoretic interest because of their close connection to special values of
Hecke's L-functions associated with quadratic imaginary fields (for this one takes lattices $\La$ with
complex multiplication).

In this paper, generalizing the works \cite{P-Wz} and \cite{P-ell}, 
we show that the Eisenstein-Kronecker numbers appear as the structure constants of
some natural $A_\infty$-algebra associated with an elliptic curve. As a consequence we derive 
quadratic relations between them which allow to express all of them polynomially in terms of the few simple ones,
which are related to the usual Eisenstein series, the Weierstrass functions $\zeta$, $\wp$ and $\wp'$, and the Kronecker functions
$F(z,w)$ and $\pa_zF(z,w)$ (see Corollary B below).
We also get a new representation for Eisenstein-Kronecker numbers by certain rapidly converging series
that appear naturally when computing the relevant $A_\infty$-products (see Theorem C).

We refer to \cite{Keller-intro} for an introduction to $A_\infty$-algebras (in particular, we use the same sign convention).
For purposes of our paper the main point is the role of an $A_\infty$-structure as a ``regularized" version of Massey products on cohomology of dg-algebras (see Sec.\ \ref{perturb-sec}). 
In particular, any $\Ext$-algebra $\Ext^*(G,G)$ of a coherent sheaf $G$ on an algebraic variety
carries such an $A_\infty$-structure. There is a choice involved in constructing it, but the result is unique up to higher homotopies.

In this paper we consider two collections of (pairwise non-isomorphic) line bundles on the
complex elliptic curve $C=\C/\La$:
$(P_i)_{i=1,\ldots,r}$ and $(L_k)_{j=1,\ldots,s}$, where $\deg P_i=0$ and $\deg L_j=1$.
We set
$$G=\bigoplus_{i=1}^r P_i\oplus \bigoplus_{j=1}^s L_j$$ 
and consider the $A_\infty$-algebra structure on $\Ext^*(G,G)$ constructed 
using the natural hermitian metrics on the relevant line bundles
(see Sec.\ \ref{perturb-sec} for details).

Assume that $\La=\Z\oplus\Z\tau$, where $\tau$ is in the upper half-plane.
We denote by
$L$ the holomorphic line bundle of degree $1$ on $C$, such that
the classical theta-function 
$$\th=\th(z,\tau)=\sum_{n\in\Z}\exp(\pi i\tau n^2+2\pi i nz)$$
descends to a global section of $L$. 
The line bundles $(P_i)$, $(L_j)$ can be written in the form
$$P_i=t_{w_i}^*L\ot L^{-1}, \ \ L_j=t_{z_j}^*L,$$
for some complex numbers $(w_i)$, $(z_j)$ (unique modulo $\La$).

We fix a generator $\xi\in H^1(C,\OO)$ which is represented by $d\ov{z}$ in the Dolbeault complex.
We denote by $\eta\in H^1(C,L^{-1})$ the unique generator such that $\eta\circ\th=\xi$.

The space $\Ext^*(G,G)$ has the following natural basis:

\noindent
(i) identity elements in $\Hom(P_i,P_i)$, $\Hom(L_j,L_j)$;

\noindent
(ii) the elements $\xi_i\in \Ext^1(P_i,P_i)$ and $\xi_j\in \Ext^1(L_j,L_j)$, corresponding to $\xi$ under
the natural isomorphism with $H^1(C,\OO)$;

\noindent
(iii) $\th_{ij}:=t_{z_j-w_i}^*\th\in H^0(t_{z_j-w_i}^*L)\simeq \Hom(P_i,L_j)$;

\noindent
(iv) $\eta_{ji}:=t_{z_j-w_i}^*\eta\in H^1(t_{z_j-w_i}^*L^{-1})\simeq \Ext^1(L_j,P_i)$.

Let us set $A=A(\La)$. It is convenient to consider the following rescaling of the elements (ii)--(iv) of the above basis:
\begin{equation}\label{twisted-basis-eq}
\begin{array}{l}
\wt{\xi}_i=A^{-1}\xi_i, \ \ \wt{\xi}_j=A^{-1}\xi_j, \ \ \wt{\th}_{ij}=\exp(A^{-1}\frac{(z_j-w_i+\ov{w_i})^2}{2})\th_{ij}, 
\\ \wt{\eta}_{ji}=A^{-1}\exp(-A^{-1}\frac{(z_j-w_i+\ov{w_i})^2}{2})\eta_{ji}.
\end{array}
\end{equation}
Note also that our $A_\infty$-structure on $\Ext^*(G,G)$ is {\it cyclic} with respect to a natural pairing 
$\lan\cdot,\cdot\ran$ on $\Ext^*(G,G)$ (see Sec.\ \ref{perturb-sec}, Eq.\ \eqref{cyclic-sym}).

Now we can describe our asnwer for the $A_\infty$-structure on $\Ext^*(G,G)$.

\medskip

\noindent
{\bf Theorem A}. {\it For $a,b,c,d\ge 0$ one has 
\begin{align*} 
&
m_n((\wt{\xi}_i)^{a},\wt{\th}_{ij},(\wt{\xi}_j)^{b},\wt{\eta}_{ji'},(\wt{\xi}_{i'})^{c}, \wt{\th}_{i'j'},(\wt{\xi}_{j'})^{d})=\\
&(-1)^{{n+1\choose 2}+1}\frac{(b+d)!}{a!b!c!d!A^{a+c}}
\cdot e^*_{a+c,b+d+1}(z_{j'}-z_j,w_i-w_{i'})\cdot \wt{\th}_{ij'}.
\end{align*} 
Note that here the indices in the pairs $(i,i')$ and $(j,j')$ are not necessarily distinct. 
The remaining $m_n$ are determined by the condition that our 
$A_\infty$-structure on $\Ext^*(G,G)$ is cyclic (see Lemma \ref{cyclic-lem}).
}

The case $r=1$ in Theorem A corresponds to the $A_\infty$-algebra associated with
an elliptic curve equipped with $n$ distinct points $p_1,\ldots,p_n$ considered in \cite{LP}. 
Indeed, the corresponding bundle $G=\OO_C\oplus \bigoplus_{j=1}^s L_j$ is related to
the object $\OO_C\oplus \bigoplus_{j=1}^s \OO_{p_j}$ considered in \cite{LP} by an autoequivalence
of $D^b(\Coh C)$.
%Fixing in addition a nonzero global holomorphic $1$-form $\la$ on $C$ one gets an isomorphism of $\Ext^*(G,G)$,
%viewed as a usual algebra, with a certain fixed algebra $E_{1,n}$ (see below).

Theorem A implies that the numbers $e^*_{a,b}(z,w)$ satisfy many quadratic equations (see Proposition \ref{ainf-id-prop}).
This leads in particular to the following result.

%These equations can also be deduced from the associative Yang-Baxter equation 
%(see \cite{P-AYBE}) satisfied by the generating function
%for these numbers, which was identified with the Kronecker function, up to a simple exponential factor, 
%in \cite[Thm.\ 1.17]{Bannai} (see Section \ref{ainf-id-sec}).

\medskip

\noindent
{\bf Corollary B}. {\it For any $a,b\ge 0$ there exist polynomials $P_{a,b}$ with rational coefficients, such that
for any $z,w\not\in\La$ one has
$$\frac{e^*_{a,b}(z,w)}{A^a}=P_{a,b}\bigl(e^*_{0,1}(z,w),e^*_{0,2}(z,w),E^*_1(z),E^*_2(z),E^*_3(z),
E^*_1(w),E^*_2(w),E^*_3(w),e^*_2,e_4,e_6\bigr),$$
where $E^*_n(z)=e^*_{0,n}(z,0)$. 
}

Note that here $e^*_{0,1}(z,w)$ differs by a simple exponential factor from the Kronecker function $F(z,w)$ (see
\eqref{e01-F-eq}, \eqref{Kron-F-eq}),
$e^*_{0,2}(z,w)=-\pa_ze^*_{0,1}(z,w)$, $E^*_1(z)=Z(z)$ is the modified Weierstrass zeta-function (see \eqref{Weier-Z-eq}), 
$E^*_2(z)$ differs from $\wp(z)$ by a constant, and $E^*_3(z)=-\wp'(z)/2$.
The polynomials $P_{a,b}$ can be computed by an explicit recursion (see Corollary \ref{recursion-cor}).

For $w=0$ (resp., $w=z=0$)
the result of Corollary B, that $e^*_{a,b}(z,0)/A^a$ (resp., $e^*_{a,b}(0,0)/A^a$)
can be expressed as a polynomial with rational coefficients
in $E^*_1(z)$, $E^*_2(z)$, $E^*_3(z)$, $e^*_2$, $e_4$ and $e_6$, follows from \cite[Prop.\ 9]{CS}
%\cite[VI.4, Eq.\ (11)]{Weil}
(resp., \cite[VI.5, p.\ 45]{Weil}).

%Note that for $z,w\not\in\La$, $e^*_{0,1}(z,w)$ is essentially the Kronecker function $F(z,w)$ (see ???), ???

Now let us describe our new series for the Eisenstein-Kronecker numbers.
For $m,n\in\Z$, a lattice $\La\sub\C$ and $z,w\in C=\C/\La$ we set
\begin{equation}\label{f-mn-def-eq}
f^*_{m,n}(z,w)=f^*_{m,n}(z,w,\La)=A^{-m}\cdot\sum_{\la\in\La\setminus\{-z\}}
\frac{(\ov{\la}+\ov{z})^m}{(\la+z)^n}\exp(-A^{-1}|\la+z|^2)\lan w,\la\ran.
\end{equation}
As in the case of $e^*_{a,b}(z,w)$, 
we use the notation $f^*$ to stress the fact that one should really consider restrictions of this function
to the strata of $C\times C$ induced by $C=\{0\}\sqcup (C\setminus\{0\})$.

\medskip

\noindent
{\bf Theorem C}. {\it One has
%$$g^*_{a,b}(z,w)=(-1)^{a+b+1}\frac{b!}{A^a}\cdot \exp(A^{-1}z(w-\ov{w}))\cdot e^*_{a,b+1}(z,-w).$$
$$e^*_{a,b+1}(z,-w)=\frac{A^a}{b!}\cdot 
\sum_{k\ge 0}k!\left((-1)^{a+b+1}{a\choose k}
f^*_{a+b-k,k+1}(w,z)\cdot\lan z,w\ran_\La+{b\choose k}f^*_{a+b-k,k+1}(z,w)\right).$$
}

Note that in the case when all $z_j$ and $w_i$ are zero (resp., $z=w=0$) Theorems A and C were proved
in \cite{P-ell}. Also, the case of $m_3$ in Theorem A (corresponding to $a=b=c=d=0$) and the identity of Theorem C
for $a=b=0$ were established in \cite{P-Wz}. The method of computation of the general $A_\infty$-products in
this paper is similar to those in \cite{P-Wz} and \cite{P-ell}. However, the answer appears in the form of a
rapidly converging series as in the right-hand side of the formula of Theorem C. 
It took the author a while to realize that these series are related to the Eisenstein-Kronecker series.

Using Theorem C we get a simpler derivation of the fact that a certain version of the Kronecker function $F$ gives
a generating series for the Eisenstein-Kronecker numbers, proved in \cite[Thm.\ 1.17]{Bannai}
(see Proposition \ref{expansion-prop}).

It is interesting to study how the $A_\infty$-structure on $\Ext^*(G,G)$
computed in Theorem A varies when we vary the parameters $(w_i)$ and $(z_j)$ determining
the line bundles $P_i$ and $L_j$ (and keeping the lattice $\La$ fixed).
There are two phenomena that we can observe using our explicit formulas.
First, the dependence of the gauge equivalence class of the $A_\infty$-structure on $(w_i)$ and $(z_j)$ should
be holomorphic. The reason for this is that there is an alternative, purely algebraic way to construct the same
gauge equivalence class of $A_\infty$-structure using Cech resolutions (see e.g., \cite[Sec.\ 3]{P-ainf} for a similar construction).
This means that the derivatives $\pa_{\ov{w_i}}$ and $\pa_{\ov{z_j}}$ of the products $(m_\bullet)$, after an appropriate
choice of a basis in $\Ext^*(G,G)$, define the trivial class in the truncated
Hochschild cohomology $HH^2_{\le 0}$ of the $A_\infty$-algebra given
by $(m_\bullet)$. In other words we should have equations of the form
$$\pa_{\ov{w_i}}(m_\bullet)=[b_i,(m_\bullet)], \ \ \pa_{\ov{z_j}}(m_\bullet)=[b_j,(m_\bullet)],$$
for appropriate $1$-cochains $b_i$ and $b_j$ in the truncated Hochschild complex $CH^*_{\le 0}$,
where $[\cdot,\cdot]$ denotes the Gerstenhaber bracket.
We show that in our case there is a very simple choice of $(b_i)$ and $(b_j)$ (see Proposition \ref{variation-prop}).

Secondly, the fact that all the $A_\infty$-algebras in our family are derived Morita equivalent also has
an infinitesimal incarnation. Namely, the infinitesimal variation of the derived Morita equivalence class is captured
by the full Hochschild cohomology $HH^2$ of the $A_\infty$-algebra. Hence, for our family we should have equations
of the form
$$\pa_{w_i}(m_\bullet)=[c_i,(m_\bullet)], \ \ \pa_{z_j}(m_\bullet)=[c_j,(m_\bullet)],$$
for appropriate $1$-cochains $c_i$ and $c_j$ in the full Hochschild complex $CH^*$.
We show (see Proposition \ref{variation-prop}) that this is indeed the case for a very simple choice of $1$-cochains. 
%Thus, for an appropriate basis $(e_k)$ of $\Ext^*(G,G)$
%we get equations of the form
%$$\pa_{w_i}m_n(e_{k_1},\ldots,e_{k_n})=\sum_{s=0}^n (-1)^? m_{n+1}(e_{k_1},\ldots,e_{k_s},c_i,e_{k_{s+1}},\ldots,e_{k_n}),$$
%and similarly for $\pa_{z_j}$ and $c_j$.
Note that similar kind of structures appear in symplectic geometry under the name of {\it pseudo-isotopies},
when one considers the variation of the Fukaya products with respect to a variation of an almost complex structure
(see \cite{Fukaya}, \cite{Tu}).
%this kind of compatibility for families of $A_\infty$-structures was considered by Junwu Tu in \cite{Tu}.
%In a slightly different setup he showed that if such a compatibility holds then the de Rham complex on the base with values in the corresponding sheaf of $A_\infty$-algebras acquires a curved $A_\infty$-structure
%(see \cite[Thm.\ 2.6]{Tu}). Meaning???

\section{New series for the Eisenstein-Kronecker numbers}

\subsection{Basic properties of the Eisenstein-Kronecker numbers}

First, let us recall some properties of $e^*_{a,b}(z,w)$. We fix a lattice $\La$ and set $\lan\cdot,\cdot\ran=\lan\cdot,\cdot\ran_\La$,
$A=A(\La)$.

Directly from the definition (using analytic continuation) we get
$$e^*_{a,b}(-z,-w)=(-1)^{a+b}e^*_{a,b}(z,w),$$
\begin{equation}\label{e-quasiper-eq}
e^*_{a,b}(z,w+\la)=e^*_{a,b}(z,w), \ \ e^*_{a,b}(z+\la,w)=e^*_{a,b}(z,w)\cdot\lan \la,w\ran^{-1}, \ \text{ for } \la\in\La,
\end{equation}
\begin{equation}\label{e-lim-eq}
e^*_{a,b}(0,w)=\lim_{z\to 0}\left[e^*_{a,b}(z,w)-\frac{\ov{z}^a}{z^b}\right]
\end{equation}
(cf. \cite[VIII.15]{Weil}).

The functional equation for the Eisenstein-Kronecker-Lerch series (see \cite[VIII.13, Eq.\ (32)]{Weil})
implies the following identity
\begin{equation}\label{e-ab-func-eq}
(b-1)!e^*_{a,b}(z,w)=A^{a-b+1}\cdot a!e^*_{b-1,a+1}(w,z)\cdot\lan w,z\ran,
\end{equation}
for $a\ge 0$, $b\ge 1$.
We also have 
\begin{equation}\label{deriv-EK-eq}
\pa_z e^*_{a,b}(z,w)=-b e^*_{a,b+1}(z,w), \ \ \pa_w e^*_{a,b}(z,w)=-A^{-1}e^*_{a+1,b}(z,w)+A^{-1}\ov{z}e^*_{a,b}(z,w),
%\pa_{\ov{z}}e^*_{a,b}(z,w)=a e^*_{a-1,b}(z,w).
\end{equation}
where in the first formula $z\not\in\La$ and in the second formula $w\not\in\La$.
It is also known that $e^*_{0,b}(0,0)=0$ for odd $b$, 
and $e^*_{0,2k}(0,0)=e^*_{2k}$, for $k\ge 1$, where
for $k\ge 2$, $e^*_{2k}=e_{2k}$ is the Eisenstein series
$$e_{2k}(\la_1,\la_2)=\sum_{\la\in L\setminus\{0\}}\frac{1}{\la^{2k}}$$
while $e^*_{0,2}(0,0)=e^*_2$ is the modified Eisenstein series
$$e^*_2=\sum_{m\in\Z}\sum_{n;n\neq 0\text{ if }m=0}\frac{1}{(m\la_2+n\la_1)^2} -A^{-1}\frac{\bar{\la_1}}{\la_1},$$
where $\La=\Z\la_1+\Z\la_2$. More generally, for $b>a\ge 0$ one has 
\begin{equation}\label{e*-ab-def-eq}
e^*_{a,b}(0,0)=\frac{(b-a-1)!}{(b-1)!}(-\DD)^a(e^*_{b-a}),
\end{equation}
where $\DD=\ov{\la_1}\pa_{\la_1}+\ov{\la_2}\pa_{\la_2}$ (see \cite[VIII.15]{Weil}).

Some of the Eisenstein-Kronecker numbers are related to other classical functions.
For example, by \cite[VIII.2, Eq.(3), p.\ 70]{Weil}, for $\La=\Z+\Z\tau$ and $z,w\in\C\setminus \La$, we have
\begin{equation}\label{e01-F-eq}
e^*_{0,1}(z,w)=2\pi i \exp(A^{-1}z(w-\ov{w}))\cdot F(z,w,\tau),
\end{equation}
where $F(z,w,\tau)$ is the {\it Kronecker function}: 
\begin{equation}\label{Kron-F-eq}
F(z,w,\tau)=-\sum_{m\ge 0,n\ge 0}\exp\bigl(2\pi i (mn\tau+mz+nw)\bigr)+\sum_{m<0,n<0}\exp\bigl(2\pi i (mn\tau+mz+nw)\bigr).
\end{equation}

Also, we have (see \cite[VIII.14]{Weil}) for $z\not\in\La$
\begin{equation}\label{e01-Z-eq}
e^*_{0,1}(z,0)=Z(z,\La),
\end{equation}
where $Z(z,\La)$ is the $\La$-periodic modification of the Weierstrass zeta-function:
\begin{equation}\label{Weier-Z-eq}
\begin{array}{l}
Z(z,\La):=\zeta(z;\la_1,\la_2)-2z_1\zeta(\frac{\la_1}{2})-2z_2\zeta(\frac{\la_2}{2}),\\
\zeta(z)=\zeta(z,\la_1,\la_2)=\frac{1}{z}+\sum_{\la\in\La\setminus\{0\}}\left(\frac{1}{z+\la}-\frac{1}{\la}+\frac{z}{\la^2}\right),
\end{array}
\end{equation}
where $\La=\Z\la_1+\Z\la_2$, $z=z_1\la_1+z_2\la_2$ with $z_1,z_2\in\R$.
%(cf. \cite[ch.~III, formula (9)]{Weil}):
%\begin{equation}\label{zeta-exp}
%\zeta(z;\la_1,\la_2)=\frac{1}{z}-\sum_{k\ge 2}e_{2k}(\la_1,\la_2)z^{2k-1}.
%\end{equation}

%Finally, the nature of discontinuity of $e^*_{0,1}(z,0)$ at $z=0$ is also known:
%one has
%$$e^*_{0,1}(0,0)=\lim_{z\to 0}[e^*_{0,1}(z,0)-

%Note that $G_{2k}(\tau)$ for $k\ge 2$ are modular but
%$G_2(\tau)$ is not. Using notation of \cite{Katz}, $2G_2(\tau)=A_2(2\pi i,2\pi i\tau)$ and the modular correction of
%$A_2$ is
%$$-\frac{S}{12}=A_2-\frac{\pi}{a(L)}\frac{\bar{\la_1}}{\la_1}.$$
%Here $S$ is modular (but not holomorphic).

%If $F(\la_1,\la_2)$ is a function of weight $k$, $f(tau)=F(2\pi i,2\pi i\tau)$ and
%$g(\tau)=WF(2\pi i,2\pi i\tau)$ then
%$$g=(\frac{1}{2\pi i}\frac{\pa}{\pa\tau}-\frac{k}{4\pi\Im(\tau)})f.$$
%In particular, if we consider the functions $W^nG_{2k}$ (resp., $W^n\wt{G}_2$) then
%the corresponding functions of $\tau$ will be polynomials in $\frac{1}{\Im(\tau})$ with coefficients
%proportional to derivatives of $G_{2k}(\tau)$.

%$$f^*_{m,n}(0,z',\La)=A^{-m}\cdot \sum_{\la\in \La\setminus\{0\}}\frac{\ov{\la}^m}{\la^n}
%\exp(-A^{-1}|\la|^2+2\pi i E_L(\la,z')),$$
%=\lim_{z\to 0}(f_{m,n}(z,z',\La)-A^{-m}\frac{\ov{z}^m}{z^n}),$$
%$$f^s_{m,n}(z,L)=(\frac{\pi}{a(L)})^m\cdot \sum_{\la}\frac{(\ov{\la}+\ov{z})^m}{(\la+z)^n}
%\exp(-\frac{\pi}{a(L)}|\la+z|^2),$$
%where $E_L(x,y)=\Im(\ov{x}y)/a(L)$.

\subsection{Rapidly converging series representation}

It will be convenient to consider the following modified versions of the functions $f^*_{m,n}(z,w)$
(see \eqref{f-mn-def-eq}):
$$\wt{f}^*_{m,n}(z,w)=A^{-m}\cdot\exp\bigl(A^{-1}z(w-\ov{w})\bigr)\sum_{\la\in\La\setminus\{-z\}}
\frac{(\ov{\la}+\ov{z})^m}{(\la+z)^n}\exp\bigl(-A^{-1}|\la+z|^2\bigr)\lan w,\la\ran.$$

Set
\begin{equation}\label{g-ab-eq}
\begin{array}{l}
\wt{g}^*_{a,b}(z,w)=\sum_{k\ge 0}k!\left({a\choose k}\wt{f}^*_{a+b-k,k+1}(w,z)+(-1)^{a+b+1}{b\choose k}\wt{f}^*_{a+b-k,k+1}(z,w)\right),\\
g^*_{a,b}(z,w)=\exp\bigl(-A^{-1}z(w-\ov{w})\bigr)\wt{g}^*_{a,b}(z,w)=\\
\sum_{k\ge 0}k!\left({a\choose k}
f^*_{a+b-k,k+1}(w,z)\cdot\lan z,w\ran+(-1)^{a+b+1}{b\choose k}f^*_{a+b-k,k+1}(z,w)\right).
\end{array}
\end{equation}
Immediately from the definition we get the identities
\begin{equation}\label{g-quasiper-eq}
g^*_{a,b}(z,w+\la)=g^*_{a,b}(z,w), \ \ g^*_{a,b}(z+\la,w)=g^*_{a,b}(z,w)\cdot\lan\la,w\ran,
\end{equation}
\begin{equation}\label{g-func-eq}
\begin{split}
\wt{g}^*_{a,b}(z,w)=(-1)^{a+b+1}\wt{g}^*_{b,a}(w,z), \\
g^*_{a,b}(z,w)=(-1)^{a+b+1} g^*_{b,a}(w,z)\cdot \lan z,w\ran. 
\end{split}
\end{equation}
Note also that 
$$f^*_{m,n}(-z,-w)=(-1)^{m+n}f^*_{m,n}(z,w), \ \ g^*_{m,n}(-z,-w)=(-1)^{m+n+1}g^*_{m,n}(z,w).$$ 

The following two identities, that are reformulations of  \cite[Thm.\ 1, Thm.\ 2]{P-Wz}, will play a key role for us.

\begin{equation}\label{zeta-identity-eq}
g^*_{0,0}(z,0)=\wt{g}^*_{0,0}(z,0)=f^*_{0,1}(0,z)-f^*_{0,1}(z,0)=-Z(z,\La),
%Z(z,L)=f^s_{0,1}(z)-f^t_{0,1}(z).
\end{equation}
\begin{equation}\label{Kron-identity-eq}
\wt{g}^*_{0,0}(z,w)=-2\pi i F(z,-w,\La),
\end{equation}
where $z,w\in\C\setminus\La$.

%In the following lemma we compute the partial derivatives of $g^*_{a,b}(z,w)$.
%Note that these are defined separately for the restriction to the strata defined by the conditions for
%one of the coordinates to be in $\La$ or in $\C\setminus\La$.

\begin{lem}\label{der-lem} 
In all of the formulas below, when computing $\pa_z$ or $\pa_{\ov{z}}$ (resp., $\pa_w$ or $\pa_{\ov{w}}$)
we assume that $z\not\in\La$ (resp., $w\not\in\La$).

(i) One has
\begin{equation}
\begin{array}{l}
\pa_z f^*_{m,n}(z,w)=-f^*_{m+1,n}(z,w)-nf^*_{m,n+1}(z,w),\\
%\pa_z f^*_{m,n}=-f^*_{m+1,n}-nf^*_{m,n+1}+A^{-1}(w-\ov{w})f_{m,n},\\
\pa_w f^*_{m,n}(z,w)=f^*_{m+1,n}(z,w)-A^{-1}\ov{z}\cdot f^*_{m,n}(z,w);\\
A\pa_{\ov{z}} f^*_{m,n}(z,w)=mf^*_{m-1,n}(z,w)-f^*_{m,n-1}(z,w), \\
A\pa_{\ov{w}} f^*_{m,n}(z,w)=-f^*_{m,n-1}(z,w)+A^{-1}z\cdot f^*_{m,n}(z,w),
\end{array}
\end{equation}
where in the last two formulas $n\ge 1$.
%$$W(f_{m,n})=f_{m+2,n}+nf_{m+1,n+1}.$$

\noindent
(ii) One has
\begin{equation}\label{g-z-deriv-eq}
\pa_z g^*_{a,b}(z,w)=g^*_{a,b+1}(z,w),
\end{equation}

\noindent
(iii)  Assume that either $w\not\in\La$, $b\ge 0$; or $w\in\La$, $b\ge 1$.
Then the function $g^*_{0,b}(z,w)$ 
 is holomorphic in $z$ varying in $\C\setminus\La$.
For $z\not\in\La$ one has
$$A\pa_{\ov{z}} g^*_{0,0}(z,w)=\de_{\La}(w),$$
where $\de_{\La}$ is the delta-function of the lattice $\La$.

\noindent 
(iv) For $a\ge 0, b\ge 0$ one has
\begin{equation}\label{g-zbar-deriv-eq}
A\pa_{\ov{z}} g^*_{a,b}(z,w)=-ag^*_{a-1,b}(z,w)+\de_{a+b,0}\de_{\La}(w).
\end{equation}

\noindent
(v) One has
\begin{equation}\label{g-w-deriv-eq}
\pa_w g^*_{a,b}(z,w)=-g^*_{a+1,b}(z,w)-A^{-1}\ov{z}g^*_{a,b}(z,w),
\end{equation}
$$A\pa_{\ov{w}} g^*_{a,b}(z,w)=bg^*_{a,b-1}(z,w)+zg^*_{a,b}(z,w)-\de_{a+b,0}\de_{\La}(z)\lan z,w\ran.$$
\end{lem}

\noindent
(vi) Let $D\sub \C$ be a small disk around $0$ such that $D\cap\La=\{0\}$. 
For $a,b\ge 0$, the function $g^*_{a,b}(z,w)+(-1)^{a+b}A^{-a}b!\cdot \frac{\ov{z}^a}{z^{b+1}}$ extends to
a real analytic function on $D\times(\C\setminus\La)$ and
\begin{equation}\label{g-lim-eq}
\left[g^*_{a,b}(z,w)+(-1)^{a+b}A^{-a}b!\cdot \frac{\ov{z}^a}{z^{b+1}}\right]|_{z=0}=g^*_{a,b}(0,w).
\end{equation}

\Pf . (i) This is a straightforward computation. Note that one has to consider separately the case when $w\in\La$
(for computing $\pa_z$ and $\pa_{\ov{z}}$) and the case $z\in\La$ (for computing $\pa_w$ and $\pa_{\ov{w}}$).

\noindent
(ii) Using the definition and part (i) we get
\begin{align*}
&\pa_z g^*_{a,b}(z,w)=\sum_{k\ge 0}
\Bigl(k!{a\choose k}f^*_{a+b-k+1,k+1}(w,z)\cdot\lan z,w\ran \\
&-(-1)^{a+b+1}k!{b\choose k}f^*_{a+b-k+1,k+1}(z,w)-(-1)^{a+b+1}(k+1)!{b\choose k}f^*_{a+b-k,k+2}(z,w)
\Bigr).
\end{align*}
Changing the summation variable in the last term to $k+1$ and using the identity ${b+1\choose k}={b\choose k}+
{b\choose k-1}$, we immediately identify the result with $g^*_{a,b+1}$.

\noindent
(iii) Assume first that $w\not\in \La$. Then the formula \eqref{Kron-identity-eq} shows that
$\wt{g}^*_{0,0}(z,w)$ is holomorphic in $z$, hence, the same is true for $g^*_{0,0}(z,w)$. By part (ii), 
$g^*_{0,b}(z,w)$ is obtained by applying $\pa_z$ to $g^*_{0,0}(z,w)$ $b$ times, hence, $g^*_{0,b}(z,w)$
is also holomorphic. 

Now assume that $w=0$ (the general $w\in\La$ is reduced to $w=0$). The formula
\eqref{zeta-identity-eq} identifies $-g^*_{0,0}(z,0)$ with $Z(z)$, which gives the required formula for $\pa_{\ov{z}}$ (using 
the Legendre period relation). Since $Z(z)$ differs
from the usual Weierstrass zeta-function $\zeta(z)$ by linear terms, we see that $g^*_{0,1}(z,0)=\pa_z g^*_{0,0}(z,0)$
is already holomorphic in $z$ (and differs from the Weierstrass $\wp$-function by a constant). 
Hence, $g^*_{0,b}(z,0)=\pa_z^{b-1}g^*_{0,1}(z,0)$ is also holomorphic.

\noindent
(iv) For $a=b=0$ the assertion follows from (iii), so we can assume that  $a+b\ge 1$. We have
\begin{align*}
&-A\pa_{\ov{z}}g^*_{a,b}(z,w)=
\sum_{k\ge 0}\Bigl(k!{a\choose k} f^*_{a+b-k,k}(w,z)\cdot\lan z,w\ran\\
&+(-1)^{a+b+1}k!{b\choose k}f^*_{a+b-k,k}-(-1)^{a+b+1}k!(a+b-k){b\choose k}f^*_{a+b-k-1,k+1}(z,w)\Bigr).
\end{align*}
Changing the summation variable in the last term to $k+1$ and using the identity
$$k!{b\choose k}-(k-1)!(a+b-k+1){b\choose k-1}=-(k-1)!a{b\choose k-1},$$ 
we can rewrite this as
\begin{equation}\label{dbar-g-ab-aux-eq}
-A\pa_{\ov{z}}g^*_{a,b}(z,w)=ag^*_{a-1,b}+f^*_{a+b,0}(w,z)\cdot\lan z,w\ran+(-1)^{a+b+1}f^*_{a+b,0}(z,w).
\end{equation}
Note that this still works for $a=0$ and shows that
$$-A\pa_{\ov{z}}g^*_{0,b}(z,w)=f^*_{b,0}(w,z)\cdot\lan z,w\ran+(-1)^{b+1}f^*_{b,0}(z,w).$$
By part (iii), this is zero for $b\ge 1$ and for $b=0$, $w\not\in\La$. Hence, for $a+b\ge 1$
we can simplify the right-hand side of
\eqref{dbar-g-ab-aux-eq} to $ag^*_{a-1,b}$.

\noindent
(v) This follows from (ii) and (iv) using \eqref{g-func-eq}.

\noindent
(vi) First, we observe that by definition of $g^*_{a,b}(z,w)$ there exists a polynomial expression 
$P_{a,b}(\ov{z},|z|^2)$ such that 
$g^*_{a,b}(z,w)+P_{a,b}(\ov{z},|z|^2)z^{-b-1}$ extends to a real analytic function on $D\times(\C\setminus\La)$, and
\begin{equation}\label{aux-g-lim-eq}
\left[g^*_{a,b}(z,w)+P_{a,b}(\ov{z},|z|^2)z^{-b-1}\right]|_{z=0}=g^*_{a,b}(0,w).
\end{equation}
Furthermore, in the case $a=b=0$ we just have $P_{0,0}=1$, i.e., $g^*_{0,0}(z,w)+z^{-1}$ is real analytic on
$D\times(\C\setminus\La)$. Note that by (ii) and (v), we have
\begin{equation}\label{g-der-formula}
g^*_{a,b}(z,w)=\pa_z^b(-\pa_w-A^{-1}\ov{z})^ag^*_{0,0}(z,w).
\end{equation}
It follows that $g^*_{a,b}(z,w)+(-1)^{a+b}A^{-a}b!\cdot \frac{\ov{z}^a}{z^{b+1}}$ is real analytic on 
$D\times(\C\setminus\La)$. This implies that the expression
$$\frac{(-1)^{a+b}A^{-a}b!\ov{z}^a-P_{a,b}(\ov{z},|z|^2)}{z^{b+1}}$$
extends to a real analytic function on $D\times(\C\setminus\La)$. Since the numerator is a polynomial in $\ov{z}$ and $z\ov{z}$,
this implies that
$$\left[\frac{(-1)^{a+b}A^{-a}b!\ov{z}^a-P_{a,b}(\ov{z},|z|^2)}{z^{b+1}}\right]|_{z=0}=0.$$
Together with \eqref{aux-g-lim-eq} this implies \eqref{g-lim-eq}.
\ed

Now we are ready to prove Theorem C which can be rewritten as
\begin{equation}\label{ThmB-eq}
(-1)^{a+b+1}\frac{A^a}{b!}g^*_{a,b}(z,w)=e^*_{a,b+1}(z,-w).
\end{equation}

\medskip

\noindent
{\it Proof of Theorem C}. First, we note that \eqref{ThmB-eq} holds for $z=w=0$.
Indeed, for $b\ge a\ge 0$ this follows from \cite[Cor.\ 1.2.2]{P-ell}, using 
\eqref{e*-ab-def-eq}. The case $b<a$ follows using the functional equation \eqref{e-ab-func-eq} together with \eqref{g-func-eq}.

Note also that the equation \eqref{ThmB-eq} is invariant with respect to the shifts of $z$ and $w$ by elements of
$\La$, due to \eqref{e-quasiper-eq} and \eqref{g-quasiper-eq}.

Next, we observe that \eqref{ThmB-eq} holds for $a=b=0$ and $z,w\not\in\La$. Indeed, both sides
are equal to $2\pi i \exp(-A^{-1}z(w-\ov{w}))\cdot F(z,-w,\tau)$, by \eqref{e01-F-eq} and \eqref{Kron-identity-eq}.

Hence, still for $z,w\not\in\La$, by \eqref{g-der-formula} and \eqref{deriv-EK-eq}, we have for any $a,b\ge 0$,
\begin{align*}
&g^*_{a,b}(z,w)=\pa_z^b(-\pa_w-A^{-1}\ov{z})^ag^*_{0,0}(z,w)=-\pa_z^b(-\pa_w-A^{-1}\ov{z})^a e^*_{0,1}(z,-w)=\\
&\frac{(-1)^{a+b+1}b!}{A^a}e^*_{a,b+1}(z,-w),
\end{align*}
so we checked \eqref{ThmB-eq} in this case.

Combining \eqref{e-lim-eq} and \eqref{g-lim-eq}, we can pass to the limit $z\to 0$ (by eliminating the singular terms),
and so \eqref{ThmB-eq} still holds for $z=0$ (equivalently, for $z\in\La$) and $w\not\in\La$.
By \eqref{e-ab-func-eq} and \eqref{g-func-eq}, this implies that \eqref{ThmB-eq} also holds when $w\in\La$, $z\not\in\La$.
Thus, we have covered all the cases.
\ed

The following result is equivalent (using Theorem C) to \cite[Thm.\ 1.17]{Bannai}. 

\begin{prop}\label{expansion-prop} 
One has for all $z_0,w_0\in\C$ the following Laurent series expansion of the meromorphic function in $z,w$
at $(0,0)$:
\begin{equation}\label{g-00-expansion}
%\exp(\frac{(z+z_0)\ov{w}-w\ov{z_0}}{A})g^*_{0,0}(z_0+z,w_0-w)=\sum_{a,b\ge 0}g^*_{a,b}(z_0,w_0)\frac{w^az^b}{a!b!}.
\Theta_{z_0,w_0}(z,w)=-\sum_{a,b\ge 0}g^*_{a,b}(z_0,-w_0)\frac{w^az^b}{a!b!}+\de_{\La}(w_0)\frac{1}{w}+
\de_{\La}(z_0)\frac{\lan w_0,z_0\ran}{z},
\end{equation}
where $\de_{\La}$ is the delta-function of the lattice $\La$.
\end{prop}

\Pf . Let us set $G(z,w)=-\Theta_{z_0,w_0}(z,w)$ for brefity.
Note that 
$$G(z,w)=\exp\bigl(\frac{(z+z_0)\ov{w}-w\ov{z_0}}{A}\bigr)g^*_{0,0}(z_0+z,-w_0-w).$$
Assume first that $z_0,w_0\not\in\La$.
The fact that $G$ is holomorphic in $z,w$ near $(0,0)$ follows immediately from Lemma \ref{der-lem}.
Now using \eqref{g-der-formula} we obtain
\begin{align*}
&\pa_z^b\pa_w^aG|_{z=w=0}=
\exp\bigl(\frac{(z+z_0)\ov{w}-w\ov{z_0}}{A}\bigr)\bigl(\pa_z+\frac{\ov{w}}{A}\bigr)^b
\bigl(-\pa_w-\frac{\ov{z_0}}{A}\bigr)^ag^*_{0,0}(z_0+z,-w_0-w)|_{z=w=0}\\
&=g^*_{a,b}(z_0,-w_0),
\end{align*}
which gives the required expansion. 

Let us now consider the case $z_0=0$, $w_0\not\in\La$.
In this case using the relation between $g^*_{0,0}$ and the Kronecker function $F$ (see \eqref{Kron-identity-eq}) we see that
$G(z,w)+z^{-1}$ is holomorphic near $(0,0)$. Hence, $\pa_z^b\pa_w^a(G+z^{-1})$ is also holomorphic near $(0,0)$, and
we just need to calculate
$\pa_z^b\pa_w^a(G+z^{-1})|_{z=w=0}$. Using \eqref{g-z-deriv-eq} and \eqref{g-w-deriv-eq} we obtain
%Assume first that $a>0$. Then $\pa_w$ kills $z^{-1}$, so we get
\begin{align*}
&\pa_z^b\pa_w^a(G+z^{-1})|_{z=w=0}=
\left[\exp\bigl(\frac{z\ov{w}}{A}\bigr)\bigl(\pa_z+\frac{\ov{w}}{A}\bigr)^b(-\pa_w)^ag^*_{0,0}(z,-w_0-w)+
\de_{a,0}\frac{(-1)^b b!}{z^{b+1}}\right]|_{z=w=0}=\\
&\left[\exp\bigl(\frac{z\ov{w}}{A}\bigr)(-1)^a\sum_{i=0}^a{a\choose i}(A^{-1}\ov{z})^ig^*_{a-i,b}(z,-w_0-w)+
\de_{a,0}\frac{(-1)^b b!}{z^{b+1}}\right]|_{z=w=0}=g^*_{a,b}(0,-w_0),
\end{align*}
where the last equality follows from \eqref{aux-g-lim-eq}.
The cases $z_0\not\in\La$, $w_0=0$ and $z_0=w_0=0$ can be considered similarly.
\ed

\section{Calculation of the $A_{\infty}$-algebra associated with an elliptic curve}

\subsection{Homological perturbation construction}\label{perturb-sec}

The homological perturbation gives an $A_{\infty}$-structure
on the cohomology of a dg-algebra $(A,d)$ equipped with a projector
$\Pi:A\to B$ onto a subspace of $\ker(d)$ and a homotopy operator $Q$
such that $\id-\Pi=dQ+Qd$. Namely, the formula for the higher products on $B$ is given by
(see \cite{KS}, \cite{Merk}) as the following sum over trees:
$$m_n(b_1,\ldots,b_n)=-\sum_{T}\eps(T)m_T(b_1,\ldots,b_n).$$
Here $T$ runs over all oriented planar rooted $3$-valent trees with $n$ leaves
marked by $b_1,\ldots,b_n$.
%left to right, and the root marked by $\Pi$ (we draw the tree in such a way that leaves are above, and
%every vertex has two edges coming from above and one from below). 
The expression
$m_T(b_1,\ldots,b_n)$ is obtained by going down from leaves to the root, applying the
multiplication in $A$ at
every vertex, applying the operator $Q$ at every inner edge, and applying $\Pi$ at the root (see \cite[sec. 6.4]{KS} for details). 
The sign $\eps(T)$ is given by
$$\eps(T)=\prod_v (-1)^{|e_1(v)|+(|e_2(v)|-1)\deg(e_1(v))},$$
where $v$ runs through vertices of $T$, $(e_1(v),e_2(v))$ is the pair of edges above $v$,
for an edge $e$ we denote by $|e|$ the total number of leaves above $e$ and by $\deg(e)$ the sum of
degrees of markings of all the leaves above $e$.

In our case the algebra
$$E=\Ext^*(G,G), \text{ where } G=\bigoplus_{i=1}^r P_i\oplus\bigoplus_{j=1}^s L_j,$$
is obtained as the cohomology of the Dolbeault dg-algebra
$$A=(\Om^{0,*}({\EE}nd(G)),\dbar).$$

To construct the homotopy operator $Q$ on $A$, as in \cite{P-hi-pr}, \cite{P-ell}, we use 
the flat metric on $C$ and on the relevant line bundles. Namely, the hermitian metric on $L$ given by
$$(f,g)=\int_C f(z)\ov{g(z)}\exp\bigl(-2\pi \frac{y^2}{\Im(\tau)}\bigr) dxdy,$$
where $z=x+iy$. To get metrics on $L_j$ we use the translation $t_{z_j}^*$. Also, we get the induced metrics
on $P_i=t_{w_i}^*L\ot L^{-1}$.

Now $A$ is the direct sum of Dolbeault
complexes of the line bundles of the form $L_j\ot P_i^{-1}$, $P_i\ot L_j^{-1}$, $L_{j'}\ot L_j^{-1}$ and $P_{i'}\ot P_i^{-1}$,
and for each such line bundle
$M$ we define the homotopy to be the operator
$$Q_M=\dbar^*G:\Om^{0,1}(M)\to \Om^{0,0}(M),$$
where $G$ is the Green operator corresponding to the Laplacian
$\dbar^*\dbar+\dbar\dbar^*$. Then $\Pi=\id-Q_M\dbar-\dbar Q_M$ is the orthogonal projection onto the 
subspace of harmonic forms in $A$.

The elements of our basis in $E=\Ext^*(G,G)$ have natural harmonic representatives
which we denote in the same way. We only need to explain what they are for elements of degree $1$:
for $\xi$ the representative is $d\ov{z}$, while for $\eta_{ij}$ we take the
translation $t^*_{z_j-w_i}$ of  
$$\eta:=\sqrt{2\Im(\tau)}\cdot\ov{\th(z,\tau)}\exp\bigl(-2\pi\frac{\Im(z)^2}{\Im(\tau)}\bigr)d\ov{z},$$ 
which is a $(0,1)$-form with values in $L^{-1}$.

We also have a natural symmetric bilinear pairing on $A=A^0\oplus A^1$ given by 
$$\lan \a,\b\ran=\frac{1}{2i\Im(\tau)}\cdot\int_C \Tr(\a\circ \b)\we dz,$$
where $\a$ and $\b$ are homogeneous elements such that $\deg(\a)+\deg(\b)=1$,
which induces a perfect pairing between $E^0$ and $E^1$.
By \cite[Thm. 1.1]{P-hi-pr}, the $A_{\infty}$-structure on $E$ satisfies
the following cyclic symmetry: 
\begin{equation}\label{cyclic-sym}
\lan m_n(\a_1,\ldots, \a_n),\a_{n+1}\ran=(-1)^{n(\deg(\a_1)+1)}\lan m_n(\a_2,\ldots,\a_{n+1}),\a_1\ran.
\end{equation}

The only nonzero double products on $E$, aside from those involving identity elements, are 
\begin{equation}\label{m2-for}
m_2(\th_{ij},\eta_{ji})=\xi_i, \ \ m_2(\eta_{ji},\th_{ij})=\xi_j.
\end{equation}

By \cite[Lem.\ 2.1.1]{P-ell}, every higher product $m_n$ containing $\id_{P_i}$ or $\id_{L_j}$ vanishes.
Together with the cyclic symmetry \eqref{cyclic-sym} this implies that the only potentially nonzero higher products are of 
the following types:

\noindent
(I) $m_n((\xi_i)^{a},\th_{ij},(\xi_j)^{b},\eta_{ji'},(\xi_{i'})^{c}, \th_{i'j'},(\xi_{j'})^{d})\in \Hom(P_i,L_{j'})$,

\noindent
(II) $m_n((\xi_j)^{a},\eta_{ji},(\xi_i)^b,\th_{ij'},(\xi_{j'})^c,\eta_{j'i'},(\xi_{i'})^d)\in\Ext^1(L_j,P_{i'})$,

\noindent
(III) $m_n((\xi_i)^{a},\th_{ij},(\xi_j)^b,\eta_{ji'},(\xi_{i'})^c,\th_{i'j'},(\xi_{j'})^d,\eta_{j',i},(\xi_i)^e)\in\Hom(P_i,P_i)$, 

\noindent
(IV) $m_n((\xi_j)^{a},\eta_{ji},(\xi_i)^b,\th_{ij'},(\xi_{j'})^c,\eta_{j'i'},(\xi_{i'})^d,\th_{i'j},(\xi_j)^e)\in\Hom(L_j,L_j)$,

\noindent
where we denote by $(\xi)^a$ the string $(\xi,\ldots,\xi)$ with $\xi$ repeated $a$ times, and the indices in the pairs $(i,i')$,
$(j,j')$ are not necessarily distinct.

Let us set
%\begin{equation}\label{Mabcd}
$$M(a,b,c,d)_{ii'}^{jj'}:=
\lan m_n((\xi_i)^{a},\th_{ij},(\xi_j)^{b},\eta_{ji'},(\xi_{i'})^{c}, \th_{i'j'},(\xi_{j'})^{d}),\eta_{j'i}\ran.
$$
%\end{equation}

\begin{lem}\label{cyclic-lem} 
The only potentially nonzero products $m_n$, with $n\ge 3$, of the $A_{\infty}$-structure on 
$E$
%B=\Ext^*(\OO\oplus \bigoplus L_i,\OO\oplus \bigoplus_i L_i)
are of the form 
$$m_n((\xi_i)^{a},\th_{ij},(\xi_j)^{b},\eta_{ji'},(\xi_{i'})^{c}, \th_{i'j'},(\xi_{j'})^{d})=M(a,b,c,d)_{ii'}^{jj'}\cdot\th_{ij'},$$
$$m_n((\xi_j)^{a},\eta_{ji},(\xi_i)^{b},\th_{ij'},(\xi_{j'})^{c}, \eta_{j'i'}, (\xi_{i'})^{d})=M(b,c,d,a)_{ii'}^{j'j}\cdot\eta_{ji'},$$
$$m_n((\xi_i)^{a},\th_{ij},(\xi_j)^{b},\eta_{ji'},(\xi_{i'})^{c}, \th_{i'j'}, (\xi_{j'})^{d},\eta_{j'i},(\xi_i)^e)=
M(a+e+1,b,c,d)_{ii'}^{jj'}\cdot\id_{P_i},$$
$$m_n((\xi_j)^{a},\eta_{ji},(\xi_i)^{b},\th_{ij'},(\xi_{j'})^{c},\eta_{j'i'},(\xi_{i'})^{d},\th_{i'j},(\xi_j)^e)=M(b,c,d,a+e+1)_{ii'}^{j'j}\cdot\id_{L_j}.$$
\end{lem}

\Pf . This follows immediately from the cyclic symmetry \eqref{cyclic-sym}.
% that gives the identities
%$$\lan m_n((\xi_j)^{a},\eta_{ji},(\xi_i)^b,\th_{ij'},(\xi_{j'})^c,\eta_{j'i'},(\xi_{i'})^d),\th_{i'j}\ran=
%\lan m_n((\xi_i)^{b},\th_{ij'},(\xi_{j'})^c,\eta_{j'i'},(\xi_{i'})^d,\th_{i'j},(\xi_{j})^a),\eta_{ji}\ran,
%$$
%\begin{align*}
%&\lan m_n((\xi_i)^{a},\th_{ij},(\xi_j)^{b},\eta_{ji'},(\xi_{i'})^{c}, \th_{i'j'}, (\xi_{j'})^{d},\eta_{j'i},(\xi_i)^e),\xi_i\ran=\\
%&\lan m_n((\xi_i)^{a+e+1},\th_{ij},(\xi_j)^b,\eta_{ji'},(\xi_{i'})^c,\th_{i'j'},(\xi_{j'})^d),\eta_{j'i}\ran,
%\end{align*}
%\begin{align*}
%&\lan m_n((\xi_j)^{a},\eta_{ji},(\xi_i)^{b},\th_{ij'},(\xi_{j'})^{c},\eta_{j'i'},(\xi_{i'})^{d},\th_{i'j},(\xi_j)^e),\xi_j\ran=\\
%&\lan m_n((\xi_i)^{b},\th_{ij'},(\xi_{j'})^c,\eta_{j'i'},(\xi_{i'})^d,\th_{i'j},(\xi_j)^{a+e+1}),\eta_{ji}\ran.
%\end{align*}
\ed

Thus, it is enough to compute the products of type (I), i.e., the
coefficients
$$\lan m_n((\xi_i)^{a},\th_{ij},(\xi_j)^b,\eta_{ji'},(\xi_{i'})^c,\th_{i'j'},(\xi_{j'})^d),\eta_{j'}\ran.$$
%$$\lan m_n((\xi)^{a},\th_i,(\xi_i)^b,\eta_i,(\xi)^c,\th_j,(\xi_j)^d),\eta_j\ran.$$

%Need to consider separately cases ???

\subsection{Calculation}

For a holomorphic line bundle $M$ with a hermitian metric we consider the operator 
$$H_M:C^{\infty}(M)\to C^{\infty}(M): s\mapsto Q_M(s\cdot d\ov{z}),$$
where $Q_M=\dbar^*G$.
Below we will denote all $\xi_i$ and $\xi_j$ simply by $\xi$.
Let us also set for $z\in \C/\La$,
$$L_z:=t_z^*L, \ \ \th_z:=t_z^*\th, \ \ \eta_z=t_z^*\eta, \ \ P_z:=t_z^*L\ot L^{-1}.$$
Note that $L_j\ot P_i^{-1}=L_{z_j-w_i}$, $P_{i'}\ot P_i^{-1}=P_{w_{i'}-w_i}$, and $L_{j'}\ot L_j^{-1}=P_{z_{j'}-z_j}$.

\begin{lem}\label{comb-lem}
One has 
\begin{align*}
&\lan m_n((\xi)^{a},\th_{ij},(\xi)^b,\eta_{ji'},(\xi)^c,\th_{i'j'},(\xi)^d),\eta_{j'i}\ran=\\
&(-1)^{{n\choose 2}+1}\cdot
\sum_{a=a_1+a_2;c=c_1+c_2}{a_2+b\choose a_2}{a_1+c_1\choose a_1}{c_2+d\choose c_2}
\Phi_{w_{i'}-w_i,z_{j'}-z_j}(a_2+b,a_1+c_1,c_2+d)+\\
&(-1)^{{n\choose 2}+n+1}\cdot\sum_{b=b_1+b_2;d=d_1+d_2}{c+d_1\choose c}{b_2+d_2\choose b_2}{a+b_1\choose a}
\Phi_{z_{j'}-z_j,w_{i'}-w_i}(c+d_1,b_2+d_2,a+b_1),
\end{align*}
where for $z_0,z'_0\in\C$ we set
$$\Phi_{z_0,z'_0}(k,l,p)=
(-1)^p\lan\Pi\left(\bigl[H_{P_{z_0}}^lQ_{P_{z_0}}(H_{L_{z_0}}^k(\th_{z_0})\cdot\eta)\bigr]\cdot H_{L_{z'_0}}^p(\th_{z'_0})\right),
\eta_{z_0+z'_0}\ran.$$
%$$\psi_{ij}(k,l,p)=(-1)^p\lan\Pi\left([H_{L_jL_i^{-1}}^lQ_{L_jL_i^{-1}}(H_{L_j}^k(\th_j)\cdot\eta_i)]\cdot H_{L_i}^p(\th_i)\right),\eta_j\ran.$$
\end{lem}

\Pf . The proof is similar to that in \cite[Sec.\ 2.3]{P-ell}.
\ed

We start by calculating the $(0,1)$-forms $H^k_{L_{z_0}}(\th_{z_0})\eta_{z'_0}$ with
values in $P_{z_0-z'_0}$. 
We use real coordinates $u,v$ on $\C$ such that $z=u+v\tau$.
Let us also set $a=\Im(\tau)$. 

Note that the function $\exp(-2\pi i z_0v)$ descends to 
an everywhere nonvanishing (but non-holomorphic) section of the line bundle $P_{z_0}$.
Hence, we have an orthonormal basis of $C^\infty$-sections of $P_{z_0}$ given by
\begin{equation}\label{Pz0-basis-eq}
\varphi_{z_0,m,n}(z)=\varphi_{z_0,m\tau-n}(z):=\exp(2\pi i (mu+nv-z_0v))=\lan m\tau-n,z\ran\cdot \exp(-2\pi i z_0v).
\end{equation}
 
\begin{lem}\label{main-lem} 
Let us set for $\la=m\tau-n$, 
$$c_\la(z)=(-1)^{mn}\exp\bigl(-\frac{\pi}{2a}(|\la|^2+2\ov{\la}z+z^2)\bigr).$$
For $k\ge 0$ and $z_0,z'_0\in\C/\La$ one has:
\begin{equation}\label{Dtheta-eq}
%H^k_{L_j}(\th_j)\eta_i=\frac{1}{k!}\cdot
H^k_{L_{z_0}}(\th_{z_0})\eta_{z'_0}=\frac{1}{k!}\cdot
\exp\bigl(2\pi i(z'_0-z_0)v'_0\bigr)\cdot
\sum_{\la\in\La}(\ov{\la}-\ov{z'_0}+\ov{z_0})^k c_{\la}(z_0-z'_0)\cdot\lan\la,z'_0\ran \cdot\varphi_{z_0-z'_0,\la}(z) d\ov{z},
\end{equation}
where $z'_0=u'_0+v'_0\tau$.
\end{lem}

\Pf . Note that 
$$H^k_{L_{z_0}}(\th_{z_0})\eta_{z'_0}=t_{z_0}^*(H^k(\th)\cdot t_{z'_0-z_0}^*\eta),$$
so our identity is equivalent to 
$$H^k(\th) t_{z'_0-z_0}^*\eta=\frac{1}{k!}\cdot
\exp\bigl(2\pi i(z'_0-z_0)(v'_0-v_0)\bigr)\cdot 
\sum_{\la}(\ov{\la}-\ov{z'_0}+\ov{z_0})^{k}c_\la(z_0-z'_0)\cdot\lan\la,z'_0-z_0\ran\cdot\varphi_{z_0-z'_0,\la}(z)
\cdot d\ov{z},$$
where $z_0=u_0+v_0\tau$.
Recall that we have 
$$H^k(\th)=\frac{1}{k!}D^k(\th),$$
where $D=-\frac{a}{\pi}\frac{\pa}{\pa z}+\ov{z}-z$ (see \cite[Lem.\ 2.4.1]{P-ell}), so our assertion is equivalent
(upon replacing $z'_0-z_0$ with $z_0$) to
\begin{equation}\label{Dtheta-bis-eq}
\begin{array}{l}
\sqrt{2a}\cdot D^k(\th)\cdot \ov{\th(z+z_0)}\exp\bigl(-2\pi a(v+v_0)^2\bigr)=\\
\exp(2\pi i z_0v_0)\cdot 
\sum_{\la}(\ov{\la}-\ov{z}_0)^{k}c_\la(-z_0)\cdot\lan\la,z_0\ran\cdot\varphi_{-z_0,\la}(z).
\end{array}
\end{equation}
For $k=0$ this follows immediately from
\cite[eq. (2.2)]{P-Wz}. To deduce the general case we use induction on $k$. Namely,
since 
$$-\frac{a}{\pi}\frac{\pa}{\pa z}\varphi_{-z_0,\la}=(\ov{\la}-z_0)\varphi_{-z_0,\la}, \ \ 
-\frac{a}{\pi}\frac{\pa}{\pa z}\exp(-2\pi av^2)=(\ov{z}-z)\exp(-2\pi av^2),$$
applying $D$ to both sides of the identity \eqref{Dtheta-bis-eq} we get
\begin{align*}
&\sqrt{2a}\cdot D^{k+1}(\th)\cdot \ov{\th(z+z_0)}\exp\bigl(-2\pi a(v+v_0)^2\bigr)+\\
&(\ov{z+z_0}-z-z_0)\cdot
\sqrt{2a}\cdot D^k(\th)\cdot \ov{\th(z+z_0)}\exp\bigl(-2\pi a(v+v_0)^2\bigr)=\exp(2\pi iz_0v_0)\times\\
&\sum_{\la}(\ov{\la}-\ov{z}_0)^{k}(\ov{\la}-z_0+\ov{z}-z)c_\la(-z_0)\cdot\lan\la,z_0\ran\cdot\varphi_{-z_0,\la}(z).
\end{align*}
It remains to use \eqref{Dtheta-bis-eq} to replace the second term in the left-hand side, and the identity for $k+1$ follows.
\ed

Now we can calculate the expressions $\Phi_{z_0,w_0}(k,l,p)$ (see Lemma \ref{comb-lem}).
Note that the calculations for $z_0\in\La$ and $z_0\not\in\La$ are slightly different, since the homotopy operator
$Q_{P_{z_0}}$ behaves differently when $P_{z_0}$ is trivial.

\begin{lem}\label{phi-lem} For either $z_0=0$ or $z_0\not\in\La$ one has 
\begin{align*}
&\Phi_{z_0,w_0}(k,l,p)=\frac{A^{l+1}}{k!p!}\cdot\exp(A^{-1} z_0(w_0-\ov{w_0}))\cdot\sum_{\la\in\La\setminus\{-z_0\}}
\frac{(\ov{\la}+\ov{z_0})^{k+p}}{(\la+z_0)^{l+1}}\exp(-A^{-1}|\la+z_0|^2)\lan -\la,w_0\ran\\
&=\frac{A^{k+p+l+1}}{k!p!}\cdot \wt{f}^*_{k+p,l+1}(z_0,w_0).
\end{align*}
%\begin{align*}
%&\Phi_{ij}(k,l,p)=\frac{1}{k!p!}(\frac{a}{\pi})^{l+1}\cdot\sum_{\la=m\tau-n\neq 0}
%\frac{\ov{\la}^{k+p}}{\la^{l+1}}\exp(-\frac{\pi}{a}|\la|^2+2\pi i[m(u_i-u_j)+n(v_i-v_j)]=\\
%&\frac{1}{k!p!}(\frac{a}{\pi})^{k+p+l+1}f^t_{k+p,l+1}(z_j-z_i),
%\end{align*}
%$$\psi_{ij}(k,l,p)=\frac{1}{k!p!}(\frac{a}{\pi})^{l+1}\cdot\sum_{\la}
%\frac{(\ov{\la}+\ov{z}_j-\ov{z}_i)^{k+p}}{(\la+z_j-z_i)^{l+1}}\exp(-\frac{\pi}{a}|\la+z_j-z_i|^2=
%\frac{1}{k!p!}(\frac{a}{\pi})^{k+p+l+1}f^s_{k+p,l+1}(z_j-z_i).$$
\end{lem}

\Pf . Assume first that $z_0\not\in\La$.
Then $P_{z_0}$ is a nontrivial line bundle of degree $0$, so that $H^0(C,L)=H^1(C,L)=0$.
Hence, in this case we have $Q_{P_{z_0}}=\dbar^{-1}$. For the elements of the basis \eqref{Pz0-basis-eq}
we have
$$\dbar \varphi_{z_0,\la}=A^{-1}(\la+z_0)\varphi_{z_0,\la}d\ov{z},$$
$$Q_{P_{z_0}}(\varphi_{z_0,\la}d\ov{z})=\frac{A}{\la+z_0}\varphi_{z_0,\la},$$
and the formula for $H_{P_{z_0}}$ is similar.

Hence, using Lemma \ref{main-lem} we obtain
$$Q_{P_{z_0}}(H_{L_{z_0}}^k(\th_{z_0})\cdot\eta)=\frac{A}{k!}\cdot
\sum_{\la}\frac{(\ov{\la}+\ov{z_0})^k}{\la+z_0} c_\la(z_0)\varphi_{z_0,\la}(z).
$$
Therefore,
\begin{equation}\label{HOQHL-eq}
%\begin{array}{l}
H_{P_{z_0}}^lQ_{P_{z_0}}(H_{L_{z_0}}^k(\th_{z_0})\cdot\eta)=
\frac{A^{l+1}}{k!}\cdot
\sum_{\la}
\frac{(\ov{\la}+\ov{z_0})^k}{(\la+z_0)^{l+1}}c_\la(z_0)\varphi_{z_0,\la}(z).
%\exp(-\frac{\pi}{2a}|m\tau-n|^2+2\pi i (mu+nv)).
%\end{array}
\end{equation}

Next, comparing the formulas for $\eta$ and for the metric on $L$ 
we observe that for a $C^{\infty}$-section $f$ of $L$ one has
$(f,\th)=0$ if and only if $\lan f,\eta\ran=0$. Hence, for $f\in\C^{\infty}(L)$ one has
$$\lan\Pi(f),\eta\ran=\lan f,\eta\ran$$
(since $\Pi$ is the orthogonal projection onto $\C\th$). The same property holds for $\eta_w\in\Om^{0,1}(L_w)$.
Therefore, 
\begin{align*}
&\Phi_{z_0,w_0}(k,l,p)=(-1)^p\lan [H^l_{P_{z_0}}Q_{P_{z_0}}(H_{L_{z_0}}^k(\th_{z_0})\cdot\eta)]\cdot H^p_{L_{w_0}}(\th_{w_0}),
\eta_{z_0+w_0}\ran=\\
&(-1)^p\lan H^l_{P_{z_0}}Q_{P_{z_0}}(H_{L_{z_0}}^k(\th_{z_0})\cdot\eta), H^p_{L_{w_0}}(\th_{w_0})\cdot\eta_{z_0+w_0}\ran.
\end{align*}
Now we observe that
$$\lan \varphi_{z_0,\la},\varphi_{-z_0,\la'}\ran=\de_{\la+\la',0}.$$
Hence, the answer can be computed using the Fourier expansion
for $H^p_{L_{w_0}}(\th_{w_0})\cdot\eta_{z_0+w_0}$ from Lemma \ref{main-lem} and the Fourier expansion \eqref{HOQHL-eq}:
\begin{align*}
&\Phi_{z_0,w_0}(k,l,p)=\frac{A^{l+1}}{k!p!}\cdot\exp(2\pi i z_0(v_0+v'_0))\sum_{\la}
\frac{(\ov{\la}+\ov{z_0})^{k+p}}{(\la+z_0)^{l+1}}c_{\la}(z_0)^2\lan -\la,z_0+w_0\ran=\\
&\frac{A^{l+1}}{k!p!}\cdot\exp(2\pi i z_0v'_0)\cdot\sum_{\la}
\frac{(\ov{\la}+\ov{z_0})^{k+p}}{(\la+z_0)^{l+1}}\exp(-A^{-1}|\la+z_0|^2)\lan -\la,w_0\ran,
\end{align*}
where $z_0=u_0+v_0\tau$, $w_0=u'_0+v'_0\tau$.

In the case $z_0\in\La$ the calculation is similar except that 
we use the fact that the operator $Q_\OO:\Om^{0,1}\to\Om^{0,0}$ is given by 
$$Q_\OO(\varphi_{0,m,n}(z)d\ov{z})=\begin{cases}\frac{a}{\pi(m\tau-n)}\varphi_{0,m,n}(z), & (m,n)\neq (0,0),\\
0 & (m,n)=(0,0).
\end{cases}
$$
\ed

\medskip

\noindent
{\it Proof of Theorem A}.
Substituting the expression for $\Phi_{z_0,w_0}(k,l,p)$
%found in Lemma \ref{phi-lem}
into the formula of Lemma \ref{comb-lem} we get
\begin{align*}
&\lan m_n((\xi)^{a},\th_{ij},(\xi)^b,\eta_{ji'},(\xi)^c,\th_{i'j'},(\xi)^d),\eta_{j'i}\ran=\\
&(-1)^{{n\choose 2}+1}A^{n-2}\cdot
\sum_{a=a_1+a_2;c=c_1+c_2}\frac{(a_1+c_1)!}{a_1!a_2!c_1!c_2!b!d!}\wt{f}^*_{a_2+c_2+b+d,a_1+c_1+1}(w_{i'}-w_i,z_{j'}-z_j)+\\
&(-1)^{{n\choose 2}+n+1}A^{n-2}\cdot
\sum_{b=b_1+b_2;d=d_1+d_2}\frac{(b_2+d_2)!}{b_1!b_2!d_1!d_2!a!c!}\wt{f}^*_{b_1+d_1+a+c,b_2+d_2+1}(z_{j'}-z_j,w_{i'}-w_i).
\end{align*}
Denoting $k=a_1+c_1$ and $l=b_2+d_2$ we can rewrite this as
\begin{align*}
&(-1)^{{n\choose 2}+1}A^{-n+2}\cdot 
\lan m_n((\xi)^{a},\th_{ij},(\xi)^b,\eta_{ji'},(\xi)^c,\th_{i'j'},(\xi)^d),\eta_{j'i}\ran=\\
&\frac{1}{b!d!}\cdot\sum_{k\ge 0}C(a,c,k)\wt{f}^*_{a+c-k+b+d,k+1}(w_{i'}-w_i,z_{j'}-z_j)+\\
&(-1)^n\frac{1}{a!c!}\cdot\sum_{l\ge 0}C(b,d,l)\wt{f}^*_{b+d-l+a+c,l+1}(z_{j'}-z_j,w_{i'}-w_i),
\end{align*}
where
$$C(a,c,k)=\sum_{a_1+c_1=k; a_1\le a, c_1\le c}\frac{k!}{a_1!c_1!(a-a_1)!(c-c_1)!}=
\frac{k!}{a!c!}\sum_{a_1+c_1=k}{a\choose a_1}{c\choose c_1}=\frac{k!}{a!c!}{a+c\choose k}$$
Thus, our formula takes form
\begin{align*}
&(-1)^{{n\choose 2}+1}A^{-n+2}\cdot 
\lan m_n((\xi)^{a},\th_{ij},(\xi)^b,\eta_{ji'},(\xi)^c,\th_{i'j'},(\xi)^d),\eta_{j'i}\ran=
\frac{1}{a!b!c!d!}\times\\
&\sum_{k\ge 0}k!\left({a+c\choose k}\wt{f}^*_{n-3-k,k+1}(w_{i'}-w_i,z_{j'}-z_j)+
(-1)^n{b+d\choose k}\wt{f}^*_{n-3-k,k+1}(z_{j'}-z_j,w_{i'}-w_i)\right).
\end{align*}
Comparing with \eqref{g-ab-eq} we obtain 
\begin{equation}\label{main-eq}
\begin{array}{l}
M(a,b,c,d)_{ii'}^{jj'}=\lan m_n((\xi)^{a},\th_{ij},(\xi)^b,\eta_{ji'},(\xi)^c,\th_{i'j'},(\xi)^d),\eta_{j'i}\ran=\\
(-1)^{{n\choose 2}+1}\frac{A^{n-2}}{a!b!c!d!}\cdot \wt{g}^*_{a+c,b+d}(z_{j'}-z_j,w_{i'}-w_i),
\end{array}
\end{equation}
where $n=a+b+c+d+3$.
Taking into account Theorem C we can rewrite this in terms of the Eisenstein-Kronecker numbers as
\begin{align*}
&m_n((\xi_i)^{a},\th_{ij},(\xi_j)^{b},\eta_{ji'},(\xi_{i'})^{c}, \th_{i'j'},(\xi_{j'})^{d})=
(-1)^{{n+1\choose 2}+1}\frac{A^{b+d+1}\cdot(b+d)!}{a!b!c!d!}\times\\
&\exp(A^{-1}(z_{j'}-z_j)(w_{i'}-w_i-\ov{w_{i'}}+\ov{w_i}))\cdot e^*_{a+c,b+d+1}(z_{j'}-z_j,w_i-w_{i'})
\cdot \th_{ij'},
\end{align*} 
which is equivalent to the formula of Theorem A.
The remaining $A_\infty$-products are determined by this, as follows from Lemma \ref{cyclic-lem}.
\ed

\subsection{Consequences of the $A_{\infty}$-constraint}\label{ainf-id-sec}

Now we will determine the identities for $g^*_{a,b}(z,w)$ (or equivalently, for $e^*_{a,b+1}(z,w)$) obtained from
the $A_{\infty}$-identities on the structure constants $M(a,b,c,d)_{ii'}^{jj'}$.

\begin{prop}\label{ainf-id-prop} 
For any integers $a,b\ge 0$ one has, for $z,z',w,w'\in\C$,
\begin{align*}
&\sum_{a=a_1+a_2}{a\choose a_1}g^*_{a_2,0}(z,w)g^*_{a_1,b}(z',w+w')-
\sum_{b=b_1+b_2}{b\choose b_1}g^*_{0,b_1}(z',w')g^*_{a,b_2}(z+z',w)+\\
&g^*_{0,0}(-z,w')g^*_{a,b}(z+z',w+w')=
\de_{b,0}\de_{\La}(w+w')\frac{1}{a+1}g^*_{a+1,0}(z,w)\\
&-\de_{a,0}\de_{\La}(z+z')\frac{\lan z+z',w\ran}{b+1}g^*_{0,b+1}(z',w')+
\de_{\La}(w)\frac{1}{a+1}g^*_{a+1,b}(z',w')\\
&-\de_{\La}(z)\lan z,w\ran g^*_{a,b+1}(z',w+w')+
\de_{\La}(w')g^*_{a+1,b}(z+z',w)-\de_{\La}(z')\frac{\lan z',w+w'\ran}{b+1}g^*_{a,b+1}(z,w).
\end{align*}
\end{prop}

\Pf .
Applying the $A_\infty$-axiom to the string 
$$(\xi_i)^a,\th_{ij},\eta_{ji'},\th_{i'j'},\eta_{j'i''},\th_{i''j''},(\xi_{j^"})^b$$
we get
\begin{align*}
&\sum_{a=a_1+a_2}(-1)^{(a_2+1)(a_1+b)+a_1}M(a_2,0,0,0)_{ii'}^{jj'}M(a_1,0,0,b)_{ii''}^{j'j''}+\\
&\sum_{b=b_1+b_2}(-1)^{(b_1+1)(a+b_2+1)+a}M(0,0,0,b_1)_{i'i''}^{j'j''}M(a,0,0,b_2)_{ii'}^{jj''}+\\
&(-1)^bM(0,0,0,0)_{i'i''}^{j'j}M(a,0,0,b)_{ii''}^{jj''}+(-1)^a\de_{b,0}\de_{i,i''}M(a+1,0,0,0)_{ii'}^{jj'}\\
&-\de_{a,0}\de_{j,j''}M(0,0,0,b+1)_{i'i''}^{j'j''}+(-1)^a[\de_{i,i'}M(a+1,0,0,b)_{i'i''}^{j'j''}
-\de_{j,j'}M(a,1,0,b)_{ii''}^{jj''}+\\
&\de_{i',i''}M(a,0,1,b)_{ii'}^{jj''}-\de_{j',j''}M(a,0,0,b+1)_{ii'}^{jj'}]
=0.
\end{align*}
Using \eqref{main-eq} we get the required identity for
$z=z_{j'}-z_j$, $z'=z_{j''}-z_{j'}$, $w=w_{i'}-w_i$, $w'=w_{i''}-w_{i'}$.
Recall that by our choice of $(w_i)$, $(z_j)$, each of these differences is either $0$ or in $\C\setminus\La$.
The general case of the identity follows using \eqref{g-quasiper-eq}.
\ed

Specializing the variables in Proposition \ref{ainf-id-prop}
by setting $w'=z'=0$, or $w=z=0$, or $w=w'=z'=0$, or $w=w'=z=0$, we get the following identities.

\begin{cor}\label{recursion-cor} 
For $z,w\not\in\La$ one has
\begin{equation}\label{g-rec1-eq}
\begin{array}{l}
g^*_{a+1,b}(z,w)-\frac{1}{b+1}g^*_{a,b+1}(z,w)=\\
\sum_{a=a_1+a_2}{a\choose a_1}g^*_{a_2,0}(z,w)g^*_{a_1,b}(0,w)-
\sum_{b=b_1+b_2}{b\choose b_1}g^*_{0,b_1}(0,0)g^*_{a,b_2}(z,w)+
g^*_{0,0}(-z,0)g^*_{a,b}(z,w);
\end{array}
\end{equation}
\begin{equation}\label{g-rec2-eq}
\begin{array}{l}
\frac{1}{a+1}g^*_{a+1,b}(z,w)-g^*_{a,b+1}(z,w)=\\
\sum_{a=a_1+a_2}{a\choose a_1}g^*_{a_2,0}(0,0)g^*_{a_1,b}(z,w)-
\sum_{b=b_1+b_2}{b\choose b_1}g^*_{0,b_1}(z,w)g^*_{a,b_2}(z,0)+
g^*_{0,0}(0,w)g^*_{a,b}(z,w);
\end{array}
\end{equation}
For $z\not\in\La$ one has
\begin{equation}\label{g-z0-rec-eq}
\begin{array}{l}
(1+\frac{\de_{b,0}}{a+1})g^*_{a+1,b}(z,0)-\frac{1}{b+1}g^*_{a,b+1}(z,0)=
\sum_{a=a_1+a_2}{a\choose a_1}g^*_{a_2,0}(z,0)g^*_{a_1,b}(0,0)\\
-\sum_{b=b_1+b_2}{b\choose b_1}g^*_{0,b_1}(0,0)g^*_{a,b_2}(z,0)+
g^*_{0,0}(-z,0)g^*_{a,b}(z,0)-\frac{1}{a+1}g^*_{a+1,b}(0,0)
\end{array}
\end{equation}
\begin{equation}\label{g-z0-rec2-eq}
\begin{array}{l}
(1+\frac{1}{a+1})g^*_{a+1,b}(z,0)-g^*_{a,b+1}(z,0)=
\sum_{a=a_1+a_2}{a\choose a_1}g^*_{a_2,0}(0,0)g^*_{a_1,b}(z,0)\\
-\sum_{b=b_1+b_2}{b\choose b_1}g^*_{0,b_1}(z,0)g^*_{a,b_2}(z,0)+
-\frac{\de_{b,0}}{a+1}g^*_{a+1,0}(0,0).
\end{array}
\end{equation}
\end{cor}

Note that for $a=0$ and $b\ge 1$ the equation \eqref{g-z0-rec2-eq} becomes
$$2g^*_{1,b}(z,0)=g^*_{0,b+1}(z,0)-\sum_{b=b_1+b_2}{b\choose b_1}g^*_{0,b_1}(z,0)g^*_{0,b_2}(z,0),$$
which is equivalent to \cite[VI.4, Eq.\ (10)]{Weil}.

\medskip

\noindent
{\it Proof of Corollary B}. The fact that all $A^{-a}e^*_{a,b}(0,0)$ can be expressed as polynomials with rational coefficients
in Eisenstein series $(e^*_n)$ is proved in \cite[VI.5]{Weil} (see also \cite[Prop.\ 2.6.1]{P-ell}). 
It is also well known that $e_n$ with $n\ge 8$ are expressed in
terms of $e_4$ and $e_6$. Thus, by \eqref{ThmB-eq}, the assertion reduces to checking that for some polynomial $P$ with
rational coefficients one has
$$g^*_{a,b}(z,w)=P\bigl(g^*_{0,0}(z,w),g^*_{0,1}(z,w),g^*_{0,0}(z),g^*_{0,1}(z),g^*_{0,2}(z),
g^*_{0,0}(w),g^*_{0,1}(w),g^*_{0,2}(w),(g^*_{m,n})),$$
where we set $g^*_{a,b}(z):=g^*_{a,b}(z,0)$, $g^*_{a,b}:=g^*_{a,b}(0,0)$. 

First, let us prove the similar assertion for $w=0$, $z\not\in\La$. In other words, we claim that there is a polynomial $P$ with
rational coefficients such that
$$g^*_{a,b}(z)=P\bigl(g^*_{0,0}(z),g^*_{0,1}(z),g^*_{0,2}(z),(g^*_{m,n})),$$
The induction on $a+b$, using identities \eqref{g-z0-rec-eq} and \eqref{g-z0-rec2-eq}, shows that we can 
express any $g^*_{a,b}(z)$ in terms of $g^*_{a',b'}(z)$ with $a'+b'\le 2$ and $(g^*_{m,n})$. 
Furthermore, using \eqref{g-z0-rec2-eq}, we can express these in terms of just $g^*_{0,0}(z)$, $g^*_{0,1}(z)$, $g^*_{0,2}(z)$
and $g^*_{m,n}$, as required.

Now in the case $z,w\not\in\La$ we similarly use \eqref{g-rec1-eq} and \eqref{g-rec2-eq} (and induction on $a+b$)
to express $g^*_{a,b}(z,w)$ in terms of $g^*_{0,0}(z,w)$, $g^*_{0,1}(z,w)$ and $(g^*_{m,n}(z,0))$.
\ed

\begin{rem} For $a=b=0$ and generic $z,w,z',w'$ the identity of Proposition \ref{ainf-id-prop} becomes
$$g^*_{0,0}(z,w)g^*_{0,0}(z',w+w')-g^*_{0,0}(z',w')g^*_{0,0}(z+z',w)+g^*_{0,0}(-z,w')g^*_{0,0}(z+z',w+w')=0$$
which means that $g^*_{0,0}$ is a scalar solution of the {\it associative Yang-Baxter equation (AYBE)}. This
is equivalent to the fact that the Kronecker function $F(z,w)$ is a solution of the AYBE (see \cite{P-AYBE}).
\end{rem}

\subsection{Variation of the $A_\infty$-structure with respect to the parameters}

Let us fix the lattice $\La$. We can view our collection of 
$A_\infty$-structures on $E=\Ext^*(G,G)$ depending on the parameters $(w_i)$, $(z_j)$,
as a bundle of minimal 
$A_\infty$-algebras on an open subset of $\C^{r+s}$ (recall that we require $w_{i'}-w_i\not\in\La$
for $i\neq i'$ and $z_{j'}-z_j\not\in\La$ for $j'\neq j$). The vector bundle is trivial with the fiber
$E$. Since we want to differentiate along the parameters on the base of the family,
the choice of a trivialization of this bundle is important for us.
We are going to use the basis \eqref{twisted-basis-eq}
in $E$ to define such a trivialization.
Note that with respect to this basis the product $m_2$ is constant in our family.

Recall that for a graded $K$-module $A$, where $K$ is a commutative ring, 
there is a graded Lie algebra on the space of Hochschild cochains given by 
the Gerstenhaber bracket. 
For $f\in\Hom_K(A^{\ot m},A)$ and $g\in \Hom_K(A^{\ot n},A)$, homogeneous of some degree
with respect to the grading on $A$, this bracket is defined as follows:
\begin{align*}
&[f,g]=f\bar{\circ} g-(-1)^{|f||g|}g\bar{\circ} f, \ \text{ where } \
f\bar{\circ}g(a_1,\ldots,a_{m+n-1})=\\
&\sum_{i=1}^m (-1)^{(|a_1|+\ldots+|a_{i-1}|+m-1)\deg(g)+(i-1)(n-1)}
f(a_1,\ldots,a_{i-1},g(a_i,\ldots,a_{i+n-1}),a_{i+n},\ldots,a_{m+n-1}),
\end{align*}
where for a cochain $f\in\Hom_K(A^{\ot m},A)$, homogeneous of degree $\deg(f)$, 
we set $|f|=\deg(f)+m-1$. Note that this bracket corresponds to the supercommutator under the
standard identification of the Hochschild cochains with coderivations of the free coalgebra cogenerated by $A[1]$.
In particular, it is a graded Lie bracket with respect to the grading $|f|$.

A collection $m=(m_n)$, where $m_n\in\Hom_K(A^{\ot n},A)$, $\deg(m_n)=2-n$, defines
an $A_\infty$-algebra over $K$ if and only if $[m,m]=0.$
The Hochschild cohomology of an $A_\infty$-algebra 
$A$ is defined as the differential of the operator $[m,?]$ on Hochschild cochains.

In our case $E$ has a natural structure of $K$-algebra for $K=\C^{r+s}$ (this corresponds to the idempotents in $E$
given by the identity morphisms of $P_i$ and $L_j$).

\begin{prop}\label{variation-prop} 
For each $j=1,\ldots,s$, $i=1,\ldots,r$, let us consider the Hochschild cochains 
$f_0(j)$, $f_0(i)$, $f_1(i)$, $f'_1(i)$, $f_2(j)$  and $f_2(i)$, where $f_a(*)\in \Hom_K(E^{\ot a},E)$ 
(depending on the parameters $(z_j)$), given by
$$f_0(j)=\wt{\xi}_j, \ \ f_0(i)=-\wt{\xi}_i,$$
$$f_1(i)(\wt{\th}_{ij})=-\ov{z_j}\wt{\th}_{ij}, \ \ f_1(i)(\wt{\eta}_{ji})=\ov{z_j}\wt{\eta}_{ji}, \ \ \
f'_1(i)(\wt{\th}_{ij})=z_j\wt{\th}_{ij}, \ \ f'_1(i)(\wt{\eta}_{ji})=-z_j\wt{\eta}_{ji},$$
$$-f_2(j)(\wt{\xi}_i,\wt{\th}_{ij})=f_2(i)(\wt{\th}_{ij},\wt{\xi}_j)=\wt{\th}_{ij}, \ \ 
f_2(j)(\wt{\eta}_{ji},\wt{\xi}_i)=-f_2(i)(\wt{\xi}_j,\wt{\eta}_{ji})=\wt{\eta}_{ji},$$
with all the other components of these cochains being zero.  
With respect to the trivialization of $E$ given by the basis \eqref{twisted-basis-eq} (along with the identity elements)
one has over the open domain where $w_{i'}-w_i\not\in\La$ for $i\neq i'$ and $z_{j'}-z_j\not\in\La$ for $j'\neq j$,
\begin{equation}\label{d-zj-m-n-eq}
\begin{array}{l}
\pa_{z_j}m_n=[m_{n+1},f_0(j)], \ \ \pa_{w_i}m_n=[m_{n+1},f_0(i)]+ [m_n,f_1(i)],\\
A\pa_{\ov{z_j}}m_n=[m_{n-1},f_2(j)], \ \ A\pa_{\ov{w_i}}m_n=[m_{n-1},f_2(i)]+ [m_n,f'_1(i)].
\end{array}
\end{equation}
Let us define a connection $\nabla$ on $E$ by requiring the basis \eqref{twisted-basis-eq} to be
horizontal with respect to $\nabla_{z_j}$ and $\nabla_{\ov{z_j}}$, and by setting
$$\nabla_{w_i}(\wt{\th}_{ij})=-A^{-1}\ov{z_j}\wt{\th}_{ij}, \ \ \nabla_{w_i}(\wt{\eta}_{ji})=A^{-1}\ov{z_j}\wt{\eta}_{ji},$$
$$\nabla_{\ov{w_i}}(\wt{\th}_{ij})=A^{-1}z_j\wt{\th}_{ij}, \ \ \nabla_{\ov{w_i}}(\wt{\eta}_{ji})=-A^{-1}z_j\wt{\eta}_{ji},$$
with other basis elements being horizontal with respect to $\nabla_{w_i}$ and $\nabla_{\ov{w_i}}$.
Then we can rewrite \eqref{d-zj-m-n-eq} as
\begin{equation}\label{d-zj-m-n-bis-eq}
\begin{array}{l}
\nabla_{z_j}m_n=[m_{n+1},f_0(j)], \ \ \nabla_{w_i}m_n=[m_{n+1},f_0(i)],\\
A\nabla_{\ov{z_j}}m_n=[m_{n-1},f_2(j)], \ \ A\nabla_{\ov{w_i}}m_n=[m_{n-1},f_2(i)].
\end{array}
\end{equation}
\end{prop}

\Pf . We use the formula of Theorem A in the form
\begin{equation}\label{aux-mn-eq}
m_n((\wt{\xi}_i)^{a},\wt{\th}_{ij},(\wt{\xi}_j)^{b},\wt{\eta}_{ji'},(\wt{\xi}_{i'})^{c}, \wt{\th}_{i'j'},(\wt{\xi}_{j'})^{d})=
\wt{M}(a,b,c,d)_{ii'}^{jj'}\cdot \wt{\th}_{ij'}, \ \text{ where}
\end{equation}
$$\wt{M}(a,b,c,d)_{ii'}^{jj'}=(-1)^{{n\choose 2}+1}\frac{1}{a!b!c!d!}\cdot g^*_{a+c,b+d}(z_{j'}-z_j,w_{i'}-w_i),$$
where $n=a+b+c+d+3$.
Now using \eqref{g-z-deriv-eq} and \eqref{g-zbar-deriv-eq} we get
$$\pa_{z_j}\wt{M}(a,b,c,d)_{ii'}^{jj'}=(-1)^{n+1}(b+1)\wt{M}(a,b+1,c,d)_{ii'}^{jj'},$$
$$\pa_{z_{j'}}\wt{M}(a,b,c,d)_{ii'}^{jj'}=(-1)^{n}(d+1)\wt{M}(a,b,c,d+1)_{ii'}^{jj'},$$
$$A\pa_{\ov{z'_j}}\wt{M}(a,b,c,d)=(-1)^n[\wt{M}(a-1,b,c,d)_{ii'}^{jj'}+\wt{M}(a,b,c-1,d)_{ii'}^{jj'}],$$
$$A\pa_{\ov{z}_j}\wt{M}(a,b,c,d)=(-1)^{n+1}[\wt{M}(a-1,b,c,d)_{ii'}^{jj'}+\wt{M}(a,b,c-1,d)_{ii'}^{jj'}],$$
$$\pa_{w_i}\wt{M}(a,b,c,d)_{ii'}^{jj'}=(-1)^{n+1}(a+1)\wt{M}(a+1,b,c,d)_{ii'}^{jj'}+A^{-1}(\ov{z_{j'}}-\ov{z_j})
\wt{M}(a,b,c,d)_{ii'}^{jj'},$$
$$\pa_{w_{i'}}\wt{M}(a,b,c,d)_{ii'}^{jj'}=(-1)^{n}(c+1)\wt{M}(a,b,c+1,d)_{ii'}^{jj'}-A^{-1}(\ov{z_{j'}}-\ov{z_j})
\wt{M}(a,b,c,d)_{ii'}^{jj'},$$
$$A\pa_{\ov{w'_i}}\wt{M}(a,b,c,d)=(-1)^{n+1}[\wt{M}(a,b-1,c,d)_{ii'}^{jj'}+\wt{M}(a,b,c,d-1)_{ii'}^{jj'}]+(z_{j'}-z_j)\wt{M}(a,b,c,d)_{ii'}^{jj'},$$
$$A\pa_{\ov{w_i}}\wt{M}(a,b,c,d)=(-1)^n[\wt{M}(a,b-1,c,d)_{ii'}^{jj'}+\wt{M}(a,b,c,d-1)_{ii'}^{jj'}]-(z_{j'}-z_j)\wt{M}(a,b,c,d)_{ii'}^{jj'}.$$
Note that in the above formulas the terms containing negative arguments in $\wt{M}(?,?,?,?)$ should be omitted.
Now the assertion follows by direct calculation (note that one should also consider other products appearing in Lemma
\ref{cyclic-lem}). For example, let us consider the equation 
$$A\nabla_{\ov{w_i}}m_n=[m_{n-1},f_2(i)]$$ 
for the products of type \eqref{aux-mn-eq}.
Since $|m_{n-1}|=2$, the right-hand side is the usual commutator, so
\begin{align*}
&[m_{n-1},f_2(i)]((\wt{\xi}_i)^{a},\wt{\th}_{ij},(\wt{\xi}_j)^{b},\wt{\eta}_{ji'},(\wt{\xi}_{i'})^{c}, \wt{\th}_{i'j'},(\wt{\xi}_{j'})^{d})=\\
&(-1)^n m_{n-1}((\wt{\xi}_i)^{a},f_2(i)(\wt{\th}_{ij},\wt{\xi}_j),(\wt{\xi}_j)^{b-1},\wt{\eta}_{ji'},(\wt{\xi}_{i'})^{c}, 
\wt{\th}_{i'j'},(\wt{\xi}_{j'})^{d})\\
&-(-1)^{n+1}f_2(i)(m_{n-1}((\wt{\xi}_i)^{a},\wt{\th}_{ij},(\wt{\xi}_j)^{b},\wt{\eta}_{ji'},(\wt{\xi}_{i'})^{c}, \wt{\th}_{i'j'},(\wt{\xi}_{j'})^{d-1}),
\wt{\xi}_{j'})=\\
&(-1)^nm_{n-1}((\wt{\xi}_i)^{a},\wt{\th}_{ij},(\wt{\xi}_j)^{b-1},\wt{\eta}_{ji'},(\wt{\xi}_{i'})^{c},\wt{\th}_{i'j'},(\wt{\xi}_{j'})^{d})+\\
&(-1)^{n}m_{n-1}((\wt{\xi}_i)^{a},\wt{\th}_{ij},(\wt{\xi}_j)^{b},\wt{\eta}_{ji'},(\wt{\xi}_{i'})^{c}, \wt{\th}_{i'j'},(\wt{\xi}_{j'})^{d-1})=\\
&(-1)^n[\wt{M}(a,b-1,c,d)_{ii'}^{jj'}+\wt{M}(a,b,c,d-1)_{ii'}^{jj'}]\wt{\th}_{ij'}.
\end{align*}
On the other hand,
\begin{align*}
&
(A\nabla_{\ov{w_i}}m_n)((\wt{\xi}_i)^{a},\wt{\th}_{ij},(\wt{\xi}_j)^{b},\wt{\eta}_{ji'},(\wt{\xi}_{i'})^{c}, \wt{\th}_{i'j'},(\wt{\xi}_{j'})^{d})=
A\pa_{\ov{w_i}}\wt{M}(a,b,c,d)_{ii'}^{jj'}\wt{\th}_{ij'}+\\
&(z_{j'}-z_j)\wt{M}(a,b,c,d)_{ii'}^{jj'}\wt{\th}_{ij'},
\end{align*}
which matches the above calculation of $[m_{n-1},f_2(i)]$, due to our formula for $A\pa_{\ov{w_i}}\wt{M}(a,b,c,d)_{ii'}^{jj'}$.
\ed

The first two equations in \eqref{d-zj-m-n-bis-eq} can be rewritted as
$$\nabla_{z_j}m_n(x_1,\ldots,x_n)=(-1)^n
\sum_{k=1}^{n+1}(-1)^{k-1+|x_1|+\ldots+|x_{k-1}|}m_{n+1}(x_1,\ldots,x_{k-1},\wt{\xi}_j,x_k,\ldots,x_n),$$
$$\nabla_{w_i}m_n(x_1,\ldots,x_n)=(-1)^{n+1}
\sum_{k=1}^{n+1}(-1)^{k-1+|x_1|+\ldots+|x_{k-1}|}m_{n+1}(x_1,\ldots,x_{k-1},\wt{\xi}_i,x_k,\ldots,x_n),$$
where $(x_i)$ are elements of the basis \eqref{twisted-basis-eq}.
Thus, in view of Theorem A we get a conceptual interpretation of the identity \eqref{g-der-formula}:
applications of $\nabla_{z_j}$ (resp., $\nabla_{w_i}$) correspond to insertions of $\wt{\xi}_j$ (resp., $\wt{\xi}_i$).
Note that equations of this form were also obtained by Tu \cite{Tu-FT} for the family of $A_\infty$-structures associated
with the fibers of a Lagrangian tori fibration.

%Should also consider $\pa_{\ov{\tau}}$???


\begin{thebibliography}{99}
\bibitem{Bannai} K.~Bannai, S.~Kobayashi, {\it Algebraic theta functions and the p-adic interpolation of Eisenstein-Kronecker numbers}, Duke Math. J. 153 (2010), no. 2, 229--295.
%\bibitem{BO} A. Bondal, D. Orlov, {\it Derived categories of coherent sheaves}, Proceedings of the International Congress of Mathematicians, Vol. II (Beijing, 2002), 47--56, Higher Ed. Press, Beijing, 
%2002.
\bibitem{CS} P.~Colmez, L.~Schneps, {\it $p$-adic interpolation of special values of Hecke L-functions},
Compositio Math. 82 (1992), 143--187.
%\bibitem{GJ} E.~Getzler, Jones
%\bibitem{Katz} N. Katz, {\it $p$-adic interpolation of real analytic
%Eisenstein series}, Annals of Math. 104 (1976), 459--571.
\bibitem{Fukaya} K.~Fukaya, {\it Cyclic symmetry and adic convergence in Lagrangian Floer theory}, Kyoto J. Math. 
50 (2010), 521--590.
\bibitem{Keller-intro} B.~Keller, {\it Introduction to $A$-infinity algebras and modules}, 
Homology Homotopy Appl. 3 (2001), no. 1, 1--35.
%\bibitem{Keller} B.~Keller, {\it $A$-infinity algebras, modules and functor categories}, in {\it Trends in representation theory of algebras and related topics}, 67--93, AMS, Providence, RI, 2006.
\bibitem{KS} M.~Kontsevich, Y.~Soibelman, {\it Homological mirror symmetry
and torus fibration}, in {\it Symplectic geometry and mirror 
symmetry (Seoul, 2000)}, 203--263, World Sci. Publishing, River Edge, NJ, 2001. 
%\bibitem{Lef}
%K.~Lef\`evre-Hasegawa, {\it Sur les $A_{\infty}$-cat\'egories}, Th\`ese de doctorat, Universit\'e Denis Diderot 
%-- Paris 7, 2003, available at B. Keller's homepage. 
\bibitem{LP} Y.~Lekili, A.~Polishchuk, 
{\it A modular compactification of ${\mathcal M}_{1,n}$ from $A_\infty$-structures}, arXiv:1408.0611.
\bibitem{Merk} S.~Merkulov, {\it Strong homotopy algebras of a K\"ahler
manifold}, IMRN 1999, no.3, 153--164.
\bibitem{P-hi-pr} A.~Polishchuk, 
         {\it Homological mirror symmetry with higher products}, 
         in {\it Proceedings of the Winter School
         on Mirror Symmetry, Vector Bundles and Lagrangian Submanifolds},
247--259, AMS and International Press, 2001. 
\bibitem{P-Wz} A.~Polishchuk, {\it Rapidly converging series for the 
Weierstrass zeta-function and for the Kronecker function}, 
 Math. Research Letters 7 (2000), 493--502.
\bibitem{P-AYBE} A.~Polishchuk, {\it Classical Yang-Baxter equation and the $A_{\infty}$-constraint},
         Advances in Math. 168 (2002), 56--95. 
\bibitem{P-ell} A.~Polishchuk, {\it $A_\infty$-algebra of an elliptic curve and Eisenstein series}, Comm. Math. Phys. 301 (2011), 709--722.
\bibitem{P-ainf} A.~Polishchuk, {\it Moduli of curves as moduli of $A_\infty$-structures}, arXiv:1312.4636.
\bibitem{Tu-FT} J.~Tu, {\it Homological mirror symmetry and Fourier-Mukai transform}, IMRN 2015, no. 3, 579--630.
\bibitem{Tu} J.~Tu, {\it On the reconstruction problem in mirror symmetry}, Adv. Math. 256 (2014), 449--478.
\bibitem{Weil} A. Weil, {\it Elliptic functions according to Eisenstein and Kronecker}. Springer-Verlag, 1976.
\end{thebibliography}
\end{document}